\documentclass[proc]{edpsmath}

\usepackage{caption}
\usepackage{subcaption}
\usepackage{amsmath,amsfonts,amsthm,amssymb}
\usepackage{mathrsfs}

\usepackage{graphicx}

\begin{document}

\title{Application of hierarchical matrix techniques to the homogenization of composite materials}
\thanks{Thanks to A. Nouy, F. Legoll, Y. Assami, A. Obliger, W. Minvielle for fruitful discussions.}
\thanks{The authors would like to CEMRACS organizers and participants for a great stay.}
\thanks{This work was financially supported by EDF R\&D.}
\author{Paul Cazeaux}\address{MATHICSE, Chair of Computational Mathematics and Simulation Sciences, Ecole Polytechnique F\'ed\'erale de Lausanne, Station 8, CH-1015 Lausanne, Szitwerland. (paul.cazeaux@epfl.ch)}
\author{Olivier Zahm}\address{GeM, UMR CNRS 6183, \'Ecole Centrale de Nantes, 1 rue de la No\"e, BP92101, 44321 Nantes Cedex 3, France. (olivier.zahm@ec-nantes.fr)}

\begin{abstract}
In this paper, we study numerical homogenization methods based on integral equations. Our work is motivated by materials such as concrete, modeled as composites structured as randomly distributed inclusions imbedded in a matrix.
We investigate two integral reformulations of the corrector problem to be solved, namely the equivalent inclusion method based on the Lippmann--Schwinger equation, and a method based on boundary integral equations.
The fully populated matrices obtained by the discretization of the integral operators are successfully dealt with using the $\mathcal{H}$-matrix format.
\end{abstract}

\begin{resume} 
Nous \'etudions la faisabilit\'e de m\'ethodes num\'eriques d'homog\'en\'eisation bas\'ees sur des \'equations int\'egrales.
Nous nous int\'eressons particuli\`erement \`a des mat\'eriaux de type b\'eton, c'est \`a dire compos\'es d'agr\'egats distribu\'es al\'eatoirement dans une matrice.
Deux reformulations \'equivalentes du probl\`eme du correcteur sont propos\'ees : la premi\`ere est la m\'ethode des inclusions \'equivalentes bas\'ee sur l'\'equation int\'egrale volumique de Lippmann--Schwinger, et la deuxi\`eme bas\'ee sur des \'equations int\'egrales surfaciques.
Les matrices pleines obtenues par la discr\'etisation des op\'erateurs int\'egraux sont trait\'ees avec succ\`es \`a l'aide de matrices hi\'erarchiques.
\end{resume}

\maketitle

\section*{Introduction}
Understanding the macroscopic properties of composite materials is of great interest in a number of industrial applications. One such example is the study of the aging of electrical nuclear plants due to the long-term behavior of concrete. Homogenization techniques have become widely used to predict this macroscopic behavior, based on the knowledge of the distribution and characteristics of the constituting elements in the microstructure of such materials. A number of different homogenization frameworks exist, ranging from qualitative methods proposed by physicists and engineering, see e.g.~\cite{Hashin1963}, to more rigorous mathematical methods developed for example in the periodic case~\cite{Bensoussan1978, Palencia1986} or the stationary random case~\cite{Bensoussan1978, RandomHomogenization}, giving rise to a vast body of literature.

Using these techniques, the effective (macroscopic) properties of the composite are typically deduced from the solution of an elliptic boundary value problem, called corrector problem. Formulated on a representative volume element (RVE) of the microstructure, this problem enables us to compute the local behavior of the microstructure under a macroscopic forcing. Approximate effective parameters are then obtained by averaging the computed local quantities over the RVE~\cite{Bensoussan1978}. In the case of idealized structures such as perfect crystals, the geometry of the RVE is simple. However, in the case of composite materials such as concrete, this geometry can become quite complex. In fact, when the composite is modeled as a matrix containing a stationary random distribution of inclusions, the corrector problem is effectively set on the whole space~\cite{RandomHomogenization}. In this case, the corrector problem is usually solved in a bounded box of increasing size. The approximate effective parameters 
thus computed are known to converge to the true effective parameters as the size of the box goes to infinity~\cite{gloria2012optimal}. Thus it is necessary to develop efficient numerical methods able to deal with complex three-dimensional geometries, induced by large numbers of inclusions in the RVE, see e.g.~\cite{Brisard11}.

A number of numerical methods can be used to approach the solution of the corrector problem. A widely used approach is the direct discretization of the elliptic PDE system, either by finite differences or by the finite elements method.
However, such full field simulations remain computationally limited by the large number of degrees of freedom that is required, in particular in a three dimensional setting.

 In this work, we propose to study two different approaches. Both are based on a reformulation of the elliptic corrector problem by integral equations. As a consequence, discretization results in linear systems involving fully populated matrices, by contrast with the sparse systems obtained with the finite differences or finite elements methods. An algebraic tool designed to deal efficiently with such fully populated matrices is the hierarchical matrix (or $\mathcal{H}$-matrix) format introduced by Hackbusch and his collaborators, see e.g.~\cite{Hackbusch1999, Borm2003, Bebendorf2008}. The main advantages of the $\mathcal{H}$-matrix format are: 
 \begin{itemize}
 \item the controlled approximation of the matrix with respect to a given precision,
 \item the reduced computational and memory cost compared to the usual matrix storage,
 \item and the accelerated algebraic operations (matrix-vector product, linear solvers, \textit{etc}).
\end{itemize} 
We investigate in this paper the applicability of this format to some linear systems, obtained by discretizing the two following integral reformulations of the corrector problem associated to the homogenization of the diffusion equation.

As a first example, we will present the equivalent inclusion method~\cite{Hashin1963,Willis1977,PONTECASTANEDA}, which is based on the Lippmann--Schwinger equation as an equivalent reformulation of the corrector problem~\cite{Dederichs1973}.
The Lippmann--Schwinger equation is a volumic integral equation with a new unknown, the polarization, set on the whole RVE. It is the starting point of analytical homogenization schemes, e.g.~\cite{Kroner1977,Willis1977,Brisard2013}, leads to the well known Hashin--Shtrikman bounds on the effective parameters~\cite{Hashin1963} and can be used to design effective numerical homogenization schemes based on the fast Fourier transform~\cite{Moulinec1994, Yvonnet2012}. The equivalent inclusion method can be seen as a crude Galerkin discretization of the Lippmann--Schwinger equation, using only constant-by-inclusion functions in the polarization space. Each entry of the corresponding matrix has an analytical expression, thus avoiding the need for numerical integration and ensuring fast assembly of the matrix~\cite{Brisard11}.

As a second example, we present a boundary integral reformulation of the corrector problem. We propose a numerical approach for its solution, using the boundary element method~\cite{TextbookBEM} and following ideas proposed by Barnett and Greengard in a different setting~\cite{BarnettGreengard}. Up to our knowledge, this work is the first attempt at solving corrector problems arising in periodic or random homogenization by use of the boundary element method. There is, however, an extensive literature on the subject of integral equations for scattering from periodic structures, see e.g.~\cite{BEMscattering}. Recently, advances have also been made in the direction of band structure calculations in periodic materials~\cite{Yuan2008, BarnettGreengard}.

The structure of the paper is as follows. In section~\ref{sec:equivalentinclusion}, we introduce the Lippmann--Schwinger equation and the equivalent inclusion method. In section~\ref{sec:hmatrix}, we recall the framework of the $\mathcal{H}$-matrix format by presenting some definitions, the main properties and associated algorithms. In section~\ref{sec:resultsequivinclusion}, we present some numerical results obtained for the equivalent inclusion method using the $\mathcal{H}$-matrix format and we assess the efficiency of this approach. Finally, in section~\ref{sec:BEM} we present the boundary integral formulation and some numerical results.

\section{Equivalent inclusion method}\label{sec:equivalentinclusion}

\subsection{Lippmann--Schwinger equation}

In this section we derive the Lippmann-Schwinger equation as an equivalent reformulation of the corrector problem associated to the homogenization of the diffusion equation~\cite{Dederichs1973}.
Let $\Omega\subset\mathbb{R}^3$ be a bounded open domain, the representative volume element of the microstructure, and $\kappa(\mathbf{x})$ a scalar diffusion coefficient. Note that we will always use the bold face convention to denote vectors.
The corrector problem is to find $u(\mathbf{x})$ such that:
\begin{equation}\label{eq:pb1}
 \left \{ \begin{aligned}
 - \mathrm{div} ( \kappa ( \nabla u+\mathbf{E} )) = 0,&&& \text{ in } \Omega,\\
 u=0 ,&&& \text{ on } \partial\Omega,
 \end{aligned} \right.
\end{equation}
where $\mathbf{E}\in\mathbb{R}^3$ is a given macroscopic potential gradient.
\begin{rmrk}
 Here we choose to work with homogeneous Dirichlet boundary conditions, but other conditions can be imposed such as periodic boundary conditions that are known to provide better estimation of the effective properties in the random homogenization framework~\cite{gloria2012optimal}.
\end{rmrk}
We introduce the polarization $\tau$ as:
\begin{equation}
 \tau(\mathbf{x}) = (\kappa(\mathbf{x})-\kappa_0)(\nabla u(\mathbf{x})+\mathbf{E}), \label{eq:deftau}
\end{equation}
where the scalar $\kappa_0>0$ corresponds to the diffusion coefficient in an homogeneous reference medium.
Then the corrector problem \eqref{eq:pb1} is equivalent to:
\begin{equation}\label{eq:pb2}
 \left \{ \begin{aligned}
 -\kappa_0 \Delta u &= \text{div}(\tau),&& \text{ in } \Omega, \\
 u&=0, &&\text{ on } \partial\Omega.
 \end{aligned} \right.
\end{equation}
Let $\delta_\mathbf{y}$ denote the Dirac mass centred in $\mathbf{y}$, and $G_0$ denote the Green's function (also called fundamental solution) associated to~\eqref{eq:pb2}, i.e. that satisfies $-\Delta G_0(\mathbf{x},\mathbf{y}) = \delta_\mathbf{y}(\mathbf{x})$ over $\Omega$, and $G_0(\mathbf{x},\mathbf{y})=0$ for $\mathbf{x}\in\partial\Omega$.
One can write:
\begin{equation}
 u(\mathbf{x})=\kappa_0^{-1}\int_\Omega G_0(\mathbf{x},\mathbf{y})\text{div}(\tau(\mathbf{y}))\text{d}\mathbf{y}=-\kappa_0^{-1}\int_\Omega \nabla_y G_0(\mathbf{x},\mathbf{y})\tau(\mathbf{x})\text{d}\mathbf{y}~~~~\forall \mathbf{x}\in\Omega.
\end{equation}
Therefore the gradient of the solution is:
\begin{equation}
 \nabla u(\mathbf{x})=-\int_\Omega \Gamma_0(\mathbf{x},\mathbf{y}):\tau(\mathbf{y})\text{d}\mathbf{y}~~~~\forall \mathbf{x}\in\Omega,
\end{equation}
where $\Gamma_0(\mathbf{x},\mathbf{y}) = \kappa_0^{-1} \nabla_\mathbf{y} \nabla_\mathbf{x} G_0(\mathbf{x},\mathbf{y})$ is a second order tensor field called the fundamental operator.
Using relation \eqref{eq:deftau}, we obtain the Lippmann-Schwinger equation for the polarization $\tau$:
\begin{equation}
 (\kappa(\mathbf{x})-\kappa_0)^{-1}:\tau(\mathbf{x}) + \int_{\Omega} \Gamma_0(\mathbf{x},\mathbf{y}):\tau(\mathbf{y})\text{d}\mathbf{y} = \mathbf{E} ~~~\forall \mathbf{x}\in \Omega. \label{eq:LS}
\end{equation}

\subsection{Energy principle of Hashin Shtrikman}
\label{sec:HashShtrik}
We look for the weak solution $\tau \in L^2(\Omega,\mathbb{R}^3)$ of \eqref{eq:LS} satisfying the variational problem $a(\tau,\widetilde\tau)=b(\widetilde\tau)$ for all $\widetilde \tau \in L^2(\Omega,\mathbb{R}^3)$, where $a(\cdot,\cdot)$ and $b(\cdot)$ are respectively the following bilinear and linear forms:
\begin{align}
 a(\tau,\widetilde \tau)&=\int_\Omega  \widetilde \tau(\mathbf{x}):(\kappa(\mathbf{x})-\kappa_0)^{-1}:\tau(\mathbf{x}) \text{d}\mathbf{x} + \int_\Omega \int_{\Omega}  \widetilde \tau(\mathbf{y}):\Gamma_0(\mathbf{x},\mathbf{y}):\tau(\mathbf{x}) \text{d}\mathbf{y}\text{d}\mathbf{x}\\
 b(\widetilde \tau)  &= \int_\Omega \widetilde \tau(\mathbf{x}):\mathbf{E}\text{d}\mathbf{x}.
\end{align}
As $\Gamma_0$ is symmetric, it was shown by Kr\"oner~\cite{Kroner1977} that this variational equation can be seen as the stationarity condition of the functional $\mathcal{H}(\tau) = \frac{1}{2}a(\tau,\tau)-b(\tau)$.
Moreover, if the relation $\kappa(\mathbf{x}) \leq \kappa_0$ (resp. $\kappa(\mathbf{x}) \geq \kappa_0$) holds in $\Omega$, the weak solution of \eqref{eq:LS} corresponds to the minimum (resp. maximum) of $\mathcal{H}$ on $L^2(\Omega,\mathbb{R}^3)$.
An approximation proposed by Willis~\cite{Willis1977} consists in replacing the fundamental operator $\Gamma_0$ by $\Gamma_\infty$, which is the fundamental operator associated to the Laplace operator in an infinite domain.
The motivation for this approximation is that $\Gamma_\infty$ has a convenient analytic expression.
The second term of $a(\tau,\tau)$ is then replaced by:
\begin{equation}
\int_\Omega \int_{\Omega}  \tau(\mathbf{y}):\Gamma_0(\mathbf{x},\mathbf{y}):\tau(\mathbf{x})  \text{d}\mathbf{y}\text{d}\mathbf{x} \approx 
\int_\Omega \int_{\Omega} \tau(\mathbf{x}):\Gamma_\infty(\mathbf{x},\mathbf{y}) : \left(\tau(\mathbf{y})-\int_{\Omega}\tau(\mathbf{z})\text{d}\mathbf{z} \right)\text{d}\mathbf{y} \text{d}\mathbf{x} . \label{eq:approxgamma}
\end{equation}
Physically this approximation holds when the effects of the boundary conditions are negligible, i.e. when the size of $\Omega$ is large. 
We obtain the new stationarity condition where $\tau^\infty$ satisfies $a^\infty(\tau^\infty,\widetilde \tau)=b(\widetilde \tau)$ for all $\widetilde \tau \in L^2(\Omega,\mathbb{R}^3)$, where
\begin{align} \label{eq:ainf}
 a^\infty(\tau^\infty,\widetilde \tau)=&\int_\Omega  \widetilde \tau(\mathbf{x}):(\kappa(\mathbf{x})-\kappa_0)^{-1}:\tau^\infty(\mathbf{x}) \text{d}\mathbf{x} 
+\int_\Omega \int_{\Omega}  \widetilde \tau(\mathbf{y}):\Gamma_\infty(\mathbf{x},\mathbf{y}):\tau^\infty(\mathbf{x}) \text{d}\mathbf{y}\text{d}\mathbf{x} \\
&- \left( \int_\Omega \widetilde\tau(\mathbf{z})\text{d}\mathbf{z} \right):\left( \int_\Omega\int_\Omega \Gamma_\infty(\mathbf{x},\mathbf{y}):\tau^\infty(\mathbf{x})\text{d}\mathbf{y} \text{d}\mathbf{x} \right). \nonumber
\end{align}
\begin{rmrk}
One can find in \cite{Brisard2013} a rigorous analysis of different approximations of $a$ corresponding to a variety of boundary conditions.
\end{rmrk}

\subsection{Assumptions on the microstructure and Galerkin approximation}

Let us now suppose that the domain $\Omega$ is spherical and is composed of $n$ spherical inclusions with diffusion coefficient $\{\kappa_\alpha\}_{\alpha\in \{1, \dots, n\}}$ imbedded in an homogeneous matrix with diffusion coefficient $\kappa_m$. Under this hypothesis, the diffusion coefficient writes:
\begin{equation}
 \kappa(\mathbf{x})= \kappa_m + \sum_{\alpha=1}^n \chi_\alpha(\mathbf{x}) (\kappa_\alpha-\kappa_m),
\end{equation}
where $\chi_\alpha(\mathbf{x})$ is the indicator function of the $\alpha^\text{th}$ inclusion.
Moreover, let us choose a reference medium with the same characteristics as the matrix, \textit{i.e. } $\kappa_0 = \kappa_m$.
By relation \eqref{eq:deftau} the polarization is zero in the matrix.
As a consequence we only need to find its value in the inclusions.

We present now the the equivalent inclusion method~\cite{Hashin1963,Willis1977,PONTECASTANEDA}. From a mathematical point of view, it consists in finding the Galerkin approximation in a finite dimensional space that is piecewise constant function on each inclusion : $V^h=\mathbb{R}^3 \otimes \text{span}(\chi_\alpha)_{1\leq \alpha\leq n}\subset L^2(\Omega,\mathbb{R}^3)$.
The coefficients $\tau_\alpha\in\mathbb{R}^3$, $\alpha\in \{1, \dots, n\}$, of the expansion $\tau^h(\mathbf{x}) = \sum_{\alpha=1}^n \tau_\alpha \chi_\alpha(\mathbf{x})$ satisfy the linear system 
\begin{equation}\label{eq:ABindiciel}
\sum_{\alpha=1}^n A_{\beta\alpha} \tau_\alpha=b_\beta~~~\forall\beta\in \{1, \dots, n\},
\end{equation}
where the 3-by-3 matrix blocs $A_{\beta\alpha}$ and the right hand side $b_\beta\in \mathbb{R}^{3}$ are:
\begin{align}
 A_{\beta\alpha} &= \delta_{\beta\alpha} f_{\beta} (\kappa_\beta - \kappa_0)^{-1} I_3
+\int_\Omega \int_{\Omega} \Gamma_\infty(\mathbf{x},\mathbf{y})\chi_{\alpha}(\mathbf{x})\chi_{\beta}(\mathbf{y})\text{d}\mathbf{x}\text{d}\mathbf{y} - f_{\beta}\int_{\Omega}\int_{\Omega}\Gamma_\infty(\mathbf{x},\mathbf{y})\chi_{\alpha}(\mathbf{x})\text{d}\mathbf{x}\text{d}\mathbf{y} , \label{eq:defA}\\
 b_{\beta} &= f_{\beta}E.
\end{align}
Here, $f_{\beta}=\int_{\Omega}\chi_\beta(\mathbf{x})\text{d}\mathbf{x}$ denotes the volume of the $\beta^{\text{th}}$ inclusion, and $I_3$ the identity matrix in $\mathbb{R}^3$.
As the domain $\Omega$ and the inclusions are spherical, the two integrals in $A_{\beta\alpha}$ have an analytic expression.
This property is an advantage of the equivalent inclusion method since there is no numerical integration to be done.
\begin{rmrk}
Note that the chosen approximation space leads to a crude approximation, since it does not take into account the variations of the polarization in the inclusions.
A more precise estimate could be reached using a more refined approximation space, for example polynomial functions of fixed degree in each inclusion, see~\cite{Moschovidis1975,Brisard11}.
\end{rmrk}
From the definition \eqref{eq:defA} it is clear that the matrix $A$ is full. Numerical resolution of the linear system \eqref{eq:ABindiciel} is then limited to a relatively small number of degrees of freedom when using classical resolution techniques.

\section{$\mathcal{H}$-matrix}\label{sec:hmatrix}

In this section we recall the main properties of the hierarchical matrix format, also called $\mathcal{H}$-matrix. This format, first introduced by Hackbusch in~\cite{Hackbusch1999}, is an algebraic tool designed mainly to manage large linear systems with fully populated matrices. For a more in-depth description, we refer to~\cite{Borm2003,Bebendorf2008}.
This format as two main advantages for numerical computation:
\begin{itemize}
 \item it provides a controlled approximation of the matrix that reduces the memory requirements,
 \item all algebraic operations (matrix-vector product, LU-factorization, \textit{etc}) can be accelerated compared to the full storage.
\end{itemize}

\subsection{Motivations}

The $\mathcal{H}$-matrix format can be introduced as a data-sparse approximation of matrices resulting from the discretization of non local integral operators~\cite{Hackbusch2002} of type:
\begin{equation}
 \int_\Omega g(\mathbf{x},\mathbf{y})u(\mathbf{y})\text{d}\mathbf{y} = f(\mathbf{x}),~~~~~\forall \mathbf{x}\in\Omega, \label{eq:integraeq}
\end{equation}
where the kernel $g$ (possibly singular) is assumed to be asymptotically smooth, that is to satisfy :
\begin{equation}
 \vert \partial_\mathbf{x}^\alpha\partial_\mathbf{y}^\beta g(\mathbf{x},\mathbf{y}) \vert \leq C_1 (C_2 \Vert \mathbf{x}-\mathbf{y} \Vert )^{-\vert\alpha\vert - \vert\beta\vert} \vert g(\mathbf{x},\mathbf{y})\vert,~~~C_1,C_2 \in \mathbb{R} ~~~~\text{with } \alpha,\beta \in \mathbb{N}^d. \label{eq:gsmooth}
\end{equation}
We look for the weak solution $u\in V$ (with $V$ an appropriate Hilbert space) of problem \eqref{eq:integraeq} which satisfies $a(u,v)=b(v)$ for all $v\in V$, where $a(\cdot,\cdot)$ and $b(\cdot)$ are defined by :
\begin{align*}
 a(u,v)&=\int_\Omega \int_\Omega g(\mathbf{x},\mathbf{y})v(\mathbf{x})u(\mathbf{y})\text{d}\mathbf{y}\text{d}\mathbf{x}, \\
 b(v)&=\int_\Omega f(\mathbf{x})v(\mathbf{x})\text{d}\mathbf{x}.
\end{align*}

Under classical assumptions on the kernel $g(\cdot,\cdot)$, the bilinear form $a$ is coercive and (by the Lax-Milgram lemma) the variational problem is well posed~\cite{Hackbusch2002}.
The Galerkin approximation $u^h$ on a finite element subspace $V^h=\text{span}\{\phi_i\}_{1\leq i\leq n}$ ($\phi_i$ being shape functions with compact, localized support) is defined by the relation $u^h = \sum_{i=1}^n U_i \phi_i$, where the coefficients $\{U_i\}_{1\leq i \leq n}$ are the solution of the linear system
\begin{equation}
AU=B, \label{eq:AUB}
\end{equation}
with $A_{i,j} = a(\phi_i,\phi_j)$ and $B_i=b(\phi_i)$.

In order to motivate the following, let us consider two subsets of indices $\tau,\sigma \subset  \{1, \dots, n\}$, and $\Omega_\tau = \cup_{i\in\tau} \text{supp}(\phi_i)$ and $\Omega_\sigma = \cup_{i\in\sigma} \text{supp}(\phi_i)$ two clusters of $\Omega$ such that :
\begin{equation}\label{eq:admissibility}
 \min \{\text{diam}(\Omega_\tau),\text{diam}(\Omega_\sigma)\} \leq \eta \text{ dist}(\Omega_\tau,\Omega_\sigma),
\end{equation}
with $\eta>0$. This last condition is called the \textit{admissibility condition}.

\begin{figure}[h]
  \centering
  \includegraphics[width=0.4\textwidth]{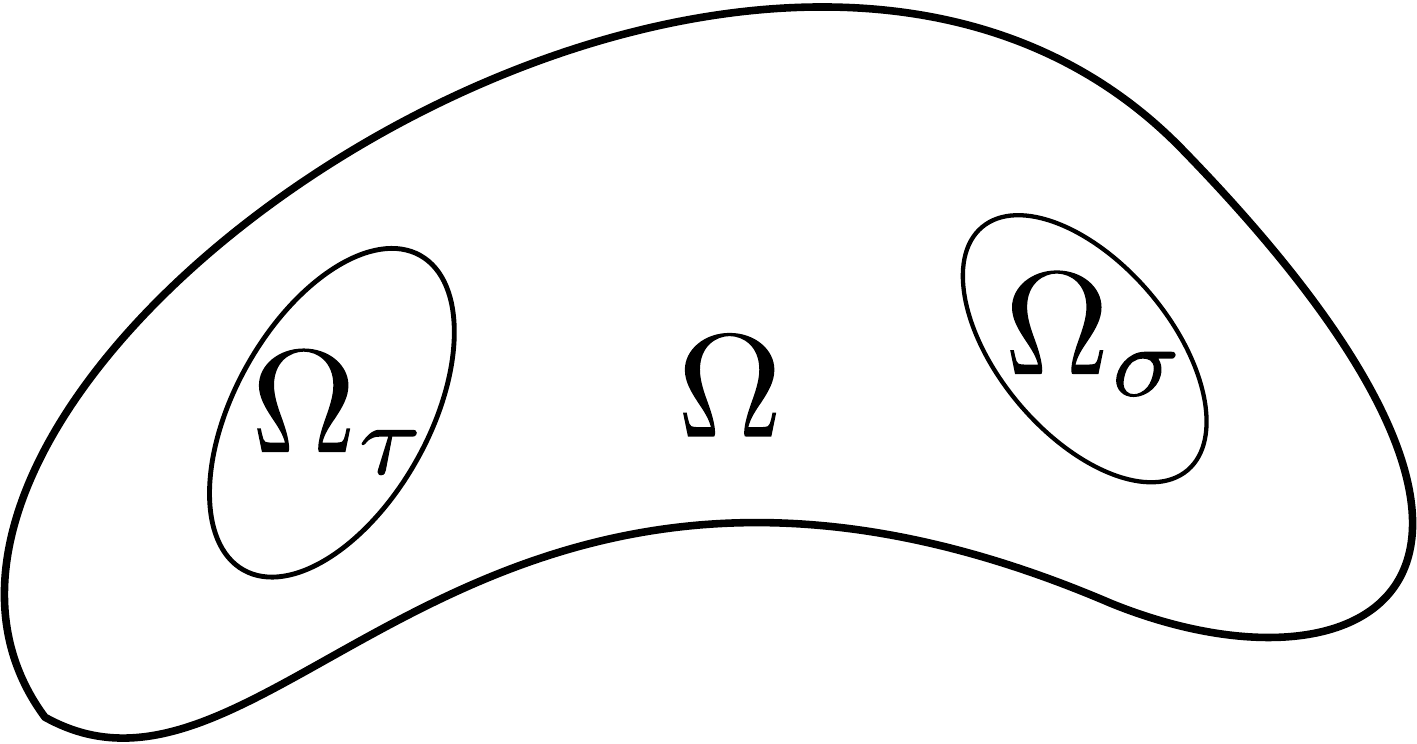}~~~~~~~
  \includegraphics[width=0.24\textwidth]{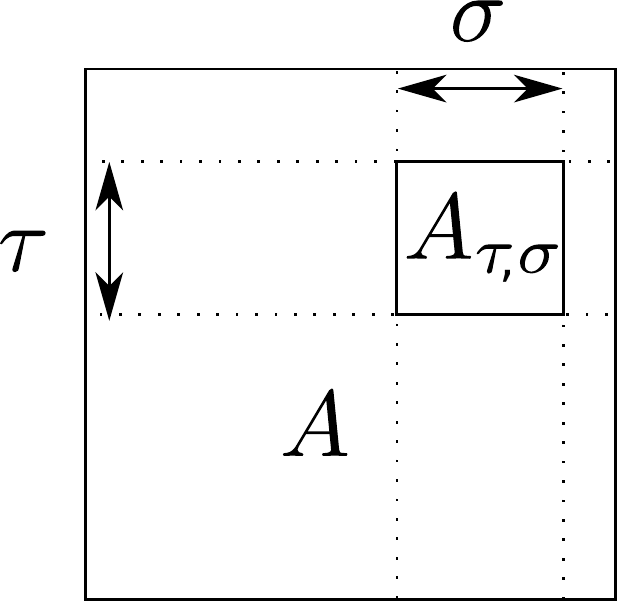}
  \caption{Representation of two admissible clusters $\Omega_\tau$ and $\Omega_\sigma$ (left) and the corresponding block in the matrix $A$ (right).}
  \label{fig:haricot}
\end{figure}

In such a situation (see figure \ref{fig:haricot}), the restriction of the kernel $g$ on $\Omega_\tau \times \Omega_\sigma$ usually has a low rank structure: for any $\varepsilon>0$ there exists $k\in\mathbb{N}$, depending on $\varepsilon$ and relatively small~\cite{Grasedyck2001}, and there exists pairs of functions $( h^1_l, h^2_l ) \in L^2(\Omega_\tau) \times L^2(\Omega_\sigma)$ with $l \in \{ 1, \dots, k \}$ such that the sum $\widetilde g(\mathbf{x},\mathbf{y}) =\sum_{l=1}^k h^1_l(\mathbf{x})h^2_l(\mathbf{y})$, called a $k$-rank function, is such that:
\begin{equation}
 \Vert g(\mathbf{x},\mathbf{y})-\widetilde g(\mathbf{x},\mathbf{y})\Vert_{L^2(\Omega_\tau \times \Omega_\sigma)} \leq \varepsilon \Vert g(\mathbf{x},\mathbf{y})\Vert_{L^2(\Omega_\tau \times \Omega_\sigma)}.
\end{equation}

For example, $\widetilde g$ can be expressed e.g. with a truncated Taylor series or an interpolation scheme, see~\cite{Borm2003,Bebendorf2008}.
The existence of such a low rank approximation ensures that the block $A_{\tau,\sigma}$ also has a low rank structure.
Indeed, replacing $g$ by $\widetilde g$ in the definition of the block $A_{\tau,\sigma}$, we obtain:
\begin{align*}
 (A_{\tau,\sigma})_{i,j} &\approx \int_\Omega \int_\Omega \widetilde g(\mathbf{x},\mathbf{y})\phi_{\sigma(j)}(\mathbf{x})\phi_{\tau(i)}(\mathbf{y})\text{d}\mathbf{y}\text{d}\mathbf{x} \\
 &= \sum_{n=1}^k \underbrace{\left( \int_\Omega h^1_n(\mathbf{x})\phi_{\sigma(j)}(\mathbf{x}) \text{d}\mathbf{x} \right)}_{\displaystyle B_{i,k}}\underbrace{\left( \int_\Omega  h^2_n(\mathbf{y})\phi_{\tau(i)}(\mathbf{y})\text{d}\mathbf{y} \right)}_{\displaystyle C_{j,k}} = (BC^t)_{i,j}.
\end{align*}
Thus, for any $\varepsilon>0$ there exists $B\in\mathbb{R}^{\vert\sigma\vert,k}$ and $C\in\mathbb{R}^{\vert\tau\vert,k}$ with $k$ small (\textit{i.e. } $k\ll \vert\sigma\vert$ and $k\ll \vert\tau\vert$) such that:
\begin{equation}
 \Vert A_{\tau,\sigma} - BC^t \Vert_F \leq \varepsilon \Vert A_{\tau,\sigma}\Vert_F, \label{eq:errorA}
\end{equation}
where $\Vert\cdot\Vert_F$ denotes the Frobenius norm, and we denote by $\vert \cdot \vert$ the cardinal of a set of indices.
The storage requirement of the approximation is $k(\vert\tau\vert+\vert\sigma\vert)$, which is significantly lower than the $\vert\tau\vert.\vert\sigma\vert$ needed for the full storage of $A_{\tau,\sigma}$.

~\\
The $\mathcal{H}$-matrix format relies on:
\begin{itemize}
 \item a block partition of the matrix that contains blocks satisfying the admissibility condition \eqref{eq:admissibility},
 \item a low-rank approximation of those admissible blocks with respect to a given precision \eqref{eq:errorA}.
\end{itemize}
In the following, we present the classical procedure for the approximation of an integral operator in the $\mathcal{H}$-matrix format~\cite{Bebendorf2008}.

\subsection{Cluster tree and block tree partition}

A good block partition of $A$ must contain a large number of admissible blocks~\cite{Bebendorf2008}. Let $I=\{1, \dots, n\}$ be the set of indices of the degrees of freedom of \eqref{eq:AUB}.
The first step is the creation of a so-called cluster tree partition $\mathcal{T}_I$ of $I$, defined in the following paragraph.
Each node of that tree is a set of indices $\sigma\subset I$ that corresponds to a subdomain $\Omega_\sigma = \cup_{i\in \sigma} \text{supp}(\phi_i)$ of $\Omega$, so that $\mathcal{T}_I$ equivalently define a partition tree of $\Omega$.
The creation of $\mathcal{T}_I$ is so that the partition of $\Omega$ contains a large number of subdomains that are potentially admissible.
The second step is the detection of admissible blocks : if $\sigma,\tau \in \mathcal{T}_I$ are so that $\Omega_\sigma$ and $\Omega_\tau$ satisfy \eqref{eq:admissibility}, then those indices corresponds to a block of $A$.
In practice, the detection of admissible blocks is done recursively in order to optimize the block partition of $A$.
It results in a block tree partition presented in the next paragraph.

\begin{rmrk}
 Only geometrical information is needed for the creation of the block partition.
\end{rmrk}

\paragraph{\bf Cluster tree partition}

A tree $\mathcal{T}_I$, with nodes $T_I$, is called a cluster tree if the following conditions hold :
\begin{enumerate}
 \item $T_I \subset \mathcal{P}(I)\backslash \{ \emptyset \}$, i.e. each node of $\mathcal{T}_I$ is a subset of the index set $I$,
 \item $I$ is the root of $\mathcal{T}_I$,
 \item If $\tau\in T_I$ is a leaf, then $\vert\tau\vert \leq C_{leaf}$, i.e. the leaves consist of a relatively small number of indices,
 \item If $\tau\in T_I$ is not a leaf, then it has two sons and their union is disjoint.
\end{enumerate}

The cluster tree $\mathcal{T}_I$ is recursively constructed with a function \textit{split}. Starting from the root $\tau = I$ and from an initial tree $\mathcal{T}_I$ that contains only the root $I$, the algorithm proceeds as follow :

\begin{enumerate}
 \item[] procedure $\mathcal{T}_I = $build\_cluster\_tree($\tau,\mathcal{T}_I$):
 \item[] if $\vert\tau\vert\geq C_{leaf}$,
 \item[] $~~~[\tau_1,\tau_2]=$ split($\tau$),
 \item[] $~~~$ add $\tau_1,$ and $\tau_2$ in $\mathcal{T}_I$ as sons of $\tau$,
 \item[] $~~~$ call $\mathcal{T}_I = $build\_cluster\_tree($\tau_1,\mathcal{T}_I$),
 \item[] $~~~$ call $\mathcal{T}_I = $build\_cluster\_tree($\tau_2,\mathcal{T}_I$),
 \item[] end.
\end{enumerate}

\begin{rmrk}
 The value $C_{leaf}=15$ classically leads to optimal computation time~\cite{Bebendorf2008}.
\end{rmrk}

\begin{figure}[h]
  \centering
  \includegraphics[scale=0.6]{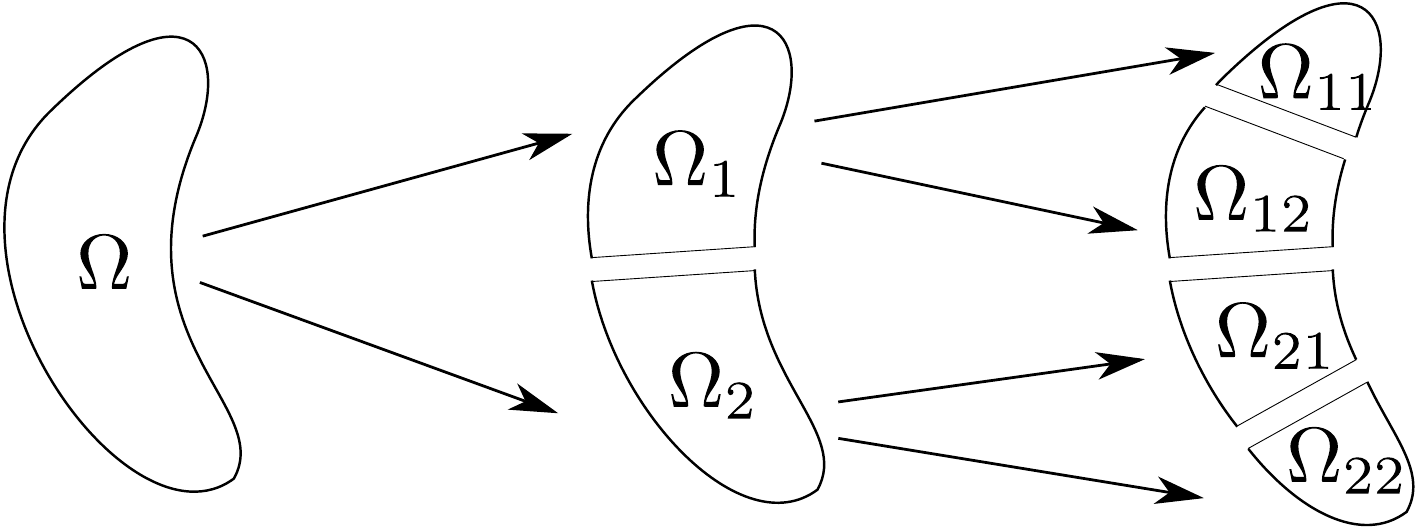}
  \caption{Construction of the cluster partitioning of $\Omega$}
  \label{fig:haricot_eclate}
\end{figure}

Note that the function $[\tau_1,\tau_2]= \text{split}(\tau)$ typically uses geometrical information to split the node $\tau$.
For each index $i\in \tau$, let $\mathbf{x}_i$ be the center of the support of $\phi_i$.
For example, a well balanced tree can be obtained with a geometric bisection on the collection of points $\mathbf{x}_i$, meaning that if $\mathbf{x}_i$ is on one side of the corresponding hyperplane, $i$ will be assigned in $\tau_1$; otherwise in $\tau_2$.
In figure \ref{fig:haricot_eclate} we illustrate the algorithm for the construction of the cluster tree partition.

\begin{figure}[h]
  \centering
  \includegraphics[scale=1]{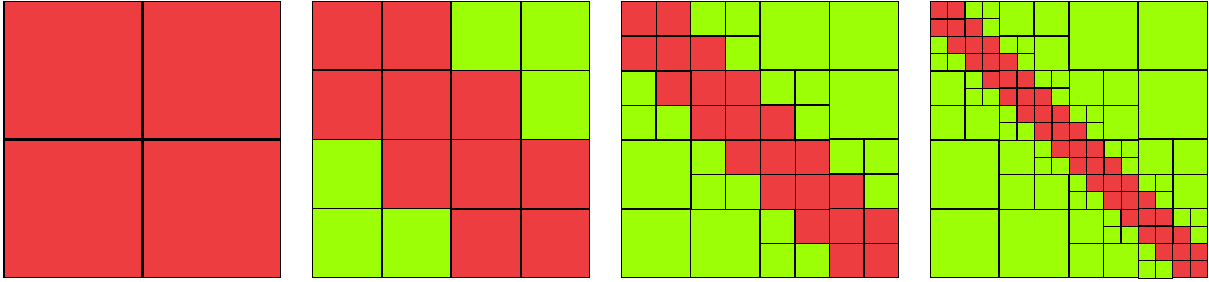}
  \caption{Construction of the block tree partition at level 1,2,3,4 of call of the \textit{build\_block\_tree} function. The green blocks are admissible : at the end, only blocks of size $<C_{leaf}$ are not admissible (pink). This matrix corresponds to a discretized 1D Laplace operator.}
  \label{fig:blocktree}
\end{figure}

\paragraph{\bf Block tree partition}
Given a cluster tree partition $\mathcal{T}_I$, the block tree partition $\mathcal{B}_{I,I}$ of $A$ is a tree that contains at each node a pair $(\tau,\sigma)$ of indices of $\mathcal{T}_I$. The leaves of $\mathcal{B}_{I,I}$ corresponds to admissible blocks \eqref{eq:admissibility}.
It can be constructed by the following procedure (see figure \ref{fig:blocktree}) initialized with $\tau=\sigma=I$ and $\mathcal{B}_{I,I}$ the block tree that contains only the root $(I,I)$:
\begin{enumerate}
 \item[] procedure $\mathcal{B}_{I,I} = $build\_block\_tree($\tau,\sigma,\mathcal{B}_{I,I}$)
 \item[] if $(\tau,\sigma)$ is not admissible, and $\vert\tau\vert \geq C_{leaf}$, and $\vert\sigma\vert \geq C_{leaf}$
 \item[] $~~~S = \{ (\tau^*,\sigma^*) , \tau^* \text{ son of }\tau, \sigma^* \text{ son of }\sigma \}$,
 \item[] $~~~$add $S$ in $\mathcal{B}_{I,I}$ as son of $(\tau,\sigma)$,
 \item[] $~~~$for $(\tau^*,\sigma^*)\in S$
 \item[] $~~~~~~\mathcal{B}_{I,I} = $build\_block\_tree($\tau^*,\sigma^*,\mathcal{B}_{I,I}$)
 \item[] $~~~$end for
 \item[] end.
\end{enumerate}

The complexity of algorithms for the creation of suitable cluster trees and block tree partitions has been analyzed in detail in~\cite{Grasedyck2001}.
For typical quasi-uniform grids, a ``good'' cluster tree can be created in $\mathcal{O}(n \log n)$ operations, the computation of the block
partition can be accomplished in $\mathcal{O}(n)$ operations.

\subsection{Approximation of the blocks}

We continue this presentation by giving here some details on the construction of a low-rank approximation of each matrix block $A_{\tau,\sigma}$.
The truncated singular value decomposition (SVD) is the optimal decomposition meaning that the relative precision $\varepsilon$ of \eqref{eq:errorA} is archieved with minimal rank $k_{\tau,\sigma}$.
But the SVD procedure requires the knowledge of the entire block $A_{\tau,\sigma}$.

To reduce the computational cost, the adaptive cross approximation (ACA) algorithm constructs a low rank approximation based on the knowledge of only a few particular rows and lines of the block. A first method for choosing these rows and lines, called ACA with full pivoting, has been proposed in~\cite{tee-cross-2000}.
In particular, it results in a quasi-optimal approximation.
The drawback of this algorithm is that the determination of the ideal pivot requires the knowledge of the entire block, as described in the next paragraph.
In~\cite{Bebendorf2000} the authors proposed a different algorithm, called ACA with partial pivoting. In this case, the choice of the pivot is made in such a way that fewer block entries are needed.
~\\
\paragraph{\bf ACA with full pivoting}

We consider a $n$-by-$m$ matrix $M$. The adaptive cross approximation is a greedy procedure on the approximation $M^k = \sum_{\nu=1}^k a^\nu\otimes b^\nu$ of $M$. Each iteration consists in the following steps : 
\begin{enumerate}
 \item Find the pivot $(i^*,j^*)$ such that :
   \begin{equation}
    (i^*,j^*)=\arg\max_{\displaystyle ij} \vert M_{ij} - M^k_{ij} \vert \label{eq:fullpivoting}
   \end{equation}
 \item Compute the two vectors $a_i^k=(M_{ij^*} - M^k_{ij^*})/(M_{i^*j^*} - M^k_{i^*j^*})$ and $b_j^k=(M_{i^*j} - M^k_{i^*j})$
 \item Update the approximation $M^{k+1}=M^k + a^k\otimes b^k$.
\end{enumerate}
The iterations are stopped when the desired precision is achieved using the condition \eqref{eq:errorA}.

In the general case there is no result on the rate of convergence.
But when the matrix corresponds to a block of an discretized integral operator with asymptotically smooth the kernel $g$, it is proved that the convergence is exponential~\cite{Bebendorf2008}. Moreover, this approximation is quasi-optimal.
However the first step is to find the largest matrix entry (full pivoting \eqref{eq:fullpivoting}), which, as for the SVD, requires the knowledge of the full matrix $M$. The ACA algorithm with partial pivoting uses a different approach for the selection of the pivot $(i^*,j^*)$ which enables a significant reduction of the computational cost.
~\\
\paragraph{\bf ACA with partial pivoting}

The idea of partial pivoting is to maximise $\vert M_{ij}-M^k_{ij}\vert$ only for one of the two indices i or j and keep the other one fixed, i.e., we determine the maximal element in modulus in one particular row or one particular column. The new pivoting strategy consists at each iteration in:
\begin{enumerate}
 \item For a given index $i^*$, find the index $j^*$ such that :
   \begin{equation}
    j^*=\arg\max_{\displaystyle j} \vert M_{i^*j}- M_{i^*j}^k \vert \label{eq:partialpivoting1}
   \end{equation}
 \item Compute the two vectors $a_i^k=(M_{ij^*} - M^k_{ij^*})/(M_{i^*j^*} - M^k_{i^*j^*})$ and $b_j^k=(M_{i^*j} - M^k_{i^*j})$
 \item Update the approximation $M^{k+1}=M^k + a^k\otimes b^k$,
 \item Find the index $i^*$ such that :
   \begin{equation}
    i^*=\arg\max_{\displaystyle i} \vert M_{ij^*}-M_{ij^*}^k \vert =\arg\max_{\displaystyle i} \vert b_j^k \vert \label{eq:partialpivoting2}
   \end{equation}
\end{enumerate}
The stopping criterion \eqref{eq:errorA} is replaced by a stagnation-based error estimator:
\begin{equation}
 \Vert a^k\otimes b^k \Vert_F \leq \varepsilon \Vert M^k\Vert_F. \label{eq:stagnation}
\end{equation}
The particularity of this algorithm is that we do not have to compute all matrix entries of $M$.
On the other hand, convergence is not guaranteed: one can find in~\cite{Borm2003} several counterexamples where ACA with partial pivoting is unable to reach the desired precision \eqref{eq:errorA}.
Note that there exists a number of variants of this algorithm, such as the improved ACA and the Hybrid Cross Approximation~\cite{Borm2003,Bebendorf2008}: using additional heuristics, they try to improve on some typical failures of the basic ACA algorithm.
However we show with numerical examples in the next section that the basic ACA algorithm seems sufficient for our application, provided that we use a well chosen value of the admissibility condition parameter~$\eta$, see~\eqref{eq:admissibility}.

\section{Numerical results}
\label{sec:resultsequivinclusion}
\begin{figure}[h]
  \centering
  \includegraphics[width=0.19\textwidth]{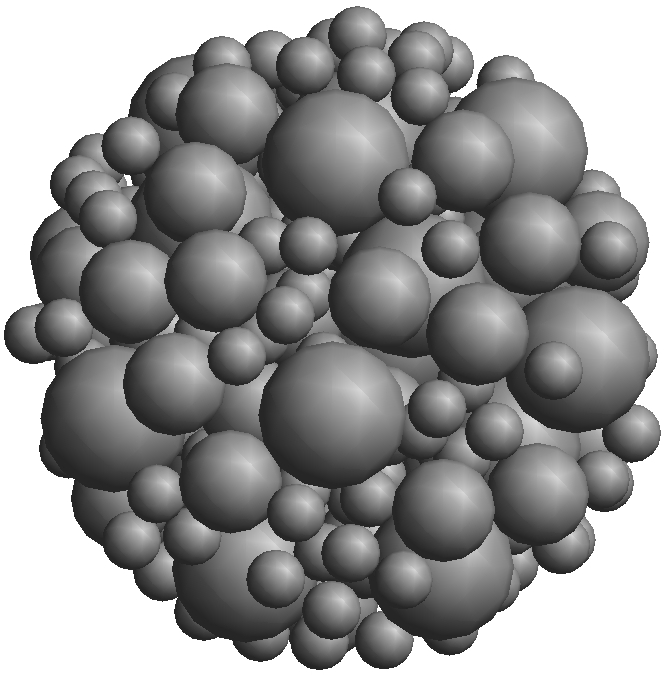}
  \includegraphics[width=0.19\textwidth]{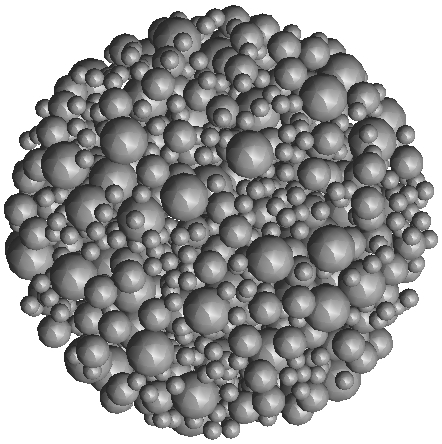}
  \includegraphics[width=0.19\textwidth]{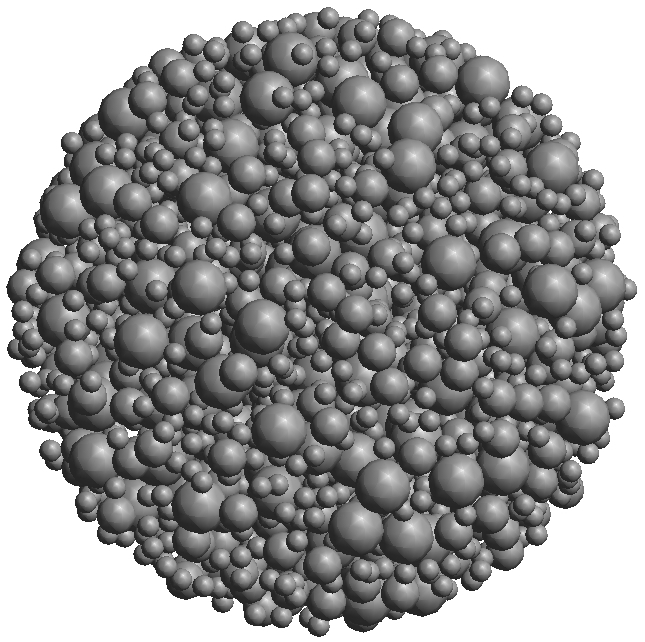}
  \includegraphics[width=0.19\textwidth]{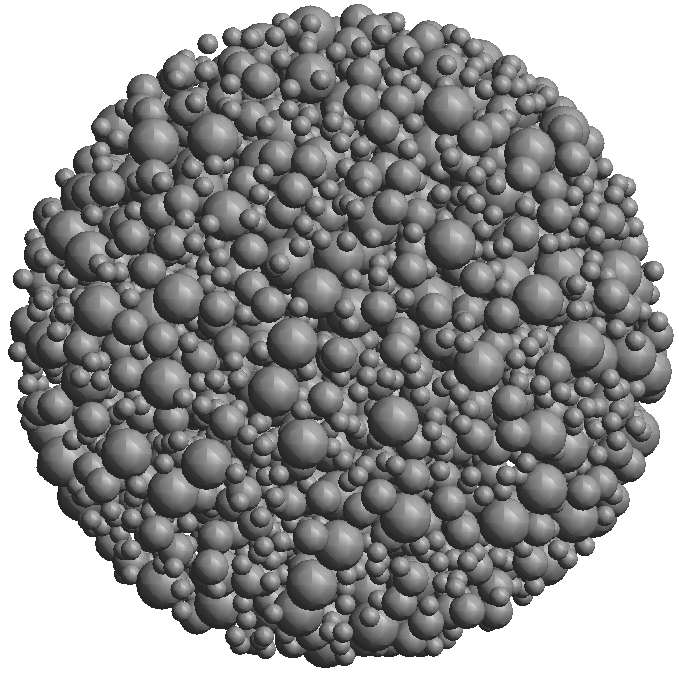}
  \includegraphics[width=0.19\textwidth]{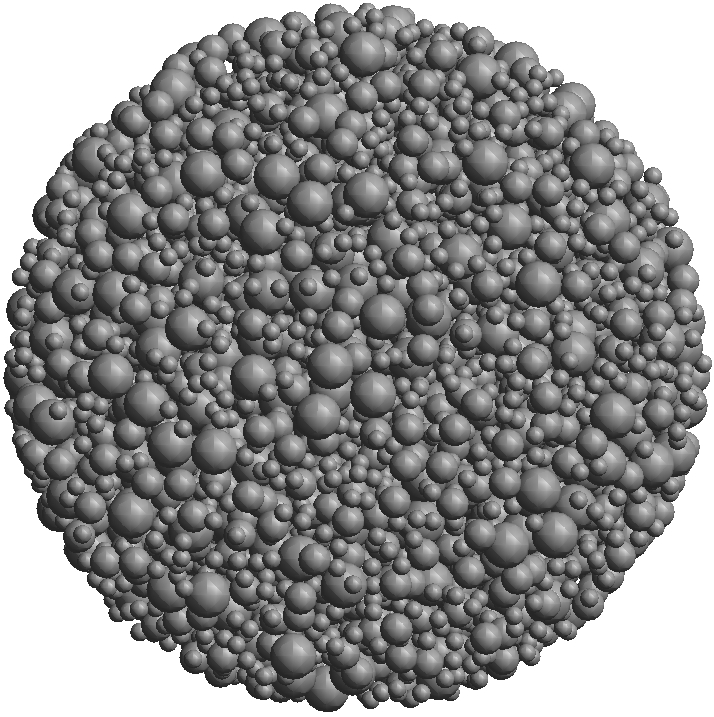}
  \caption{Representation of the five studied microstructures (200, 1000, 2000, 3000 and 5000 inclusions in a spherical domain).}
  \label{fig:four_geometries}
\end{figure}

We present in this section numerical results on the use of the $\mathcal{H}$-matrix for the resolution of the equivalent inclusion method.
The library used for the $\mathcal{H}$-matrix approximation is Ahmed~\cite{AHMED}, but others are also available such as Hlib, Hlib-pro.
This library provides all the necessary procedures to manage the $\mathcal{H}$-matrix approximation: cluster and block tree partition, ACA assembly, linear solvers, preconditioners...
Five different microstructures are studied and contain respectively 200, 1000, 2000, 3000 and 5000 spherical inclusions randomly distributed, see figure \ref{fig:four_geometries}. The coefficient diffusion of the matrix is 1, and 100 in the inclusions.
Let us denote by $A^H$ the $\mathcal{H}$-matrix approximation of the equivalent inclusion matrix $A$ \eqref{eq:ABindiciel}. Note that the ACA algorithm with partial pivoting is used for the computation of $A^H$.

\subsection{$\mathcal{H}$-matrix approximation}

\begin{figure}[h]
  \centering
  \begin{subfigure}[b]{0.245\textwidth}
   \includegraphics[width=\textwidth]{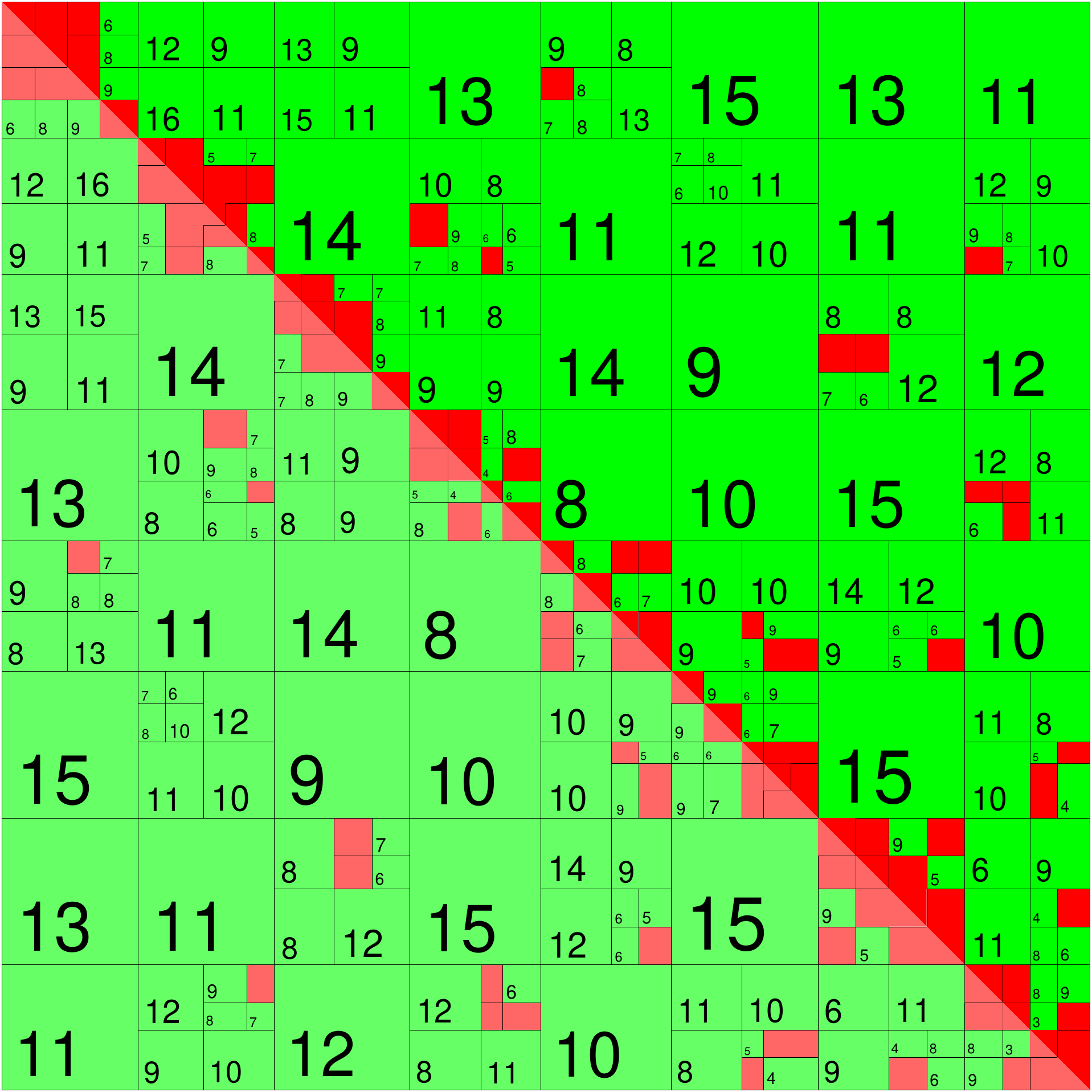}
   \caption{$\eta=1.4,\varepsilon=1e^{-2}$}
  \end{subfigure}
  \begin{subfigure}[b]{0.245\textwidth}
   \includegraphics[width=\textwidth]{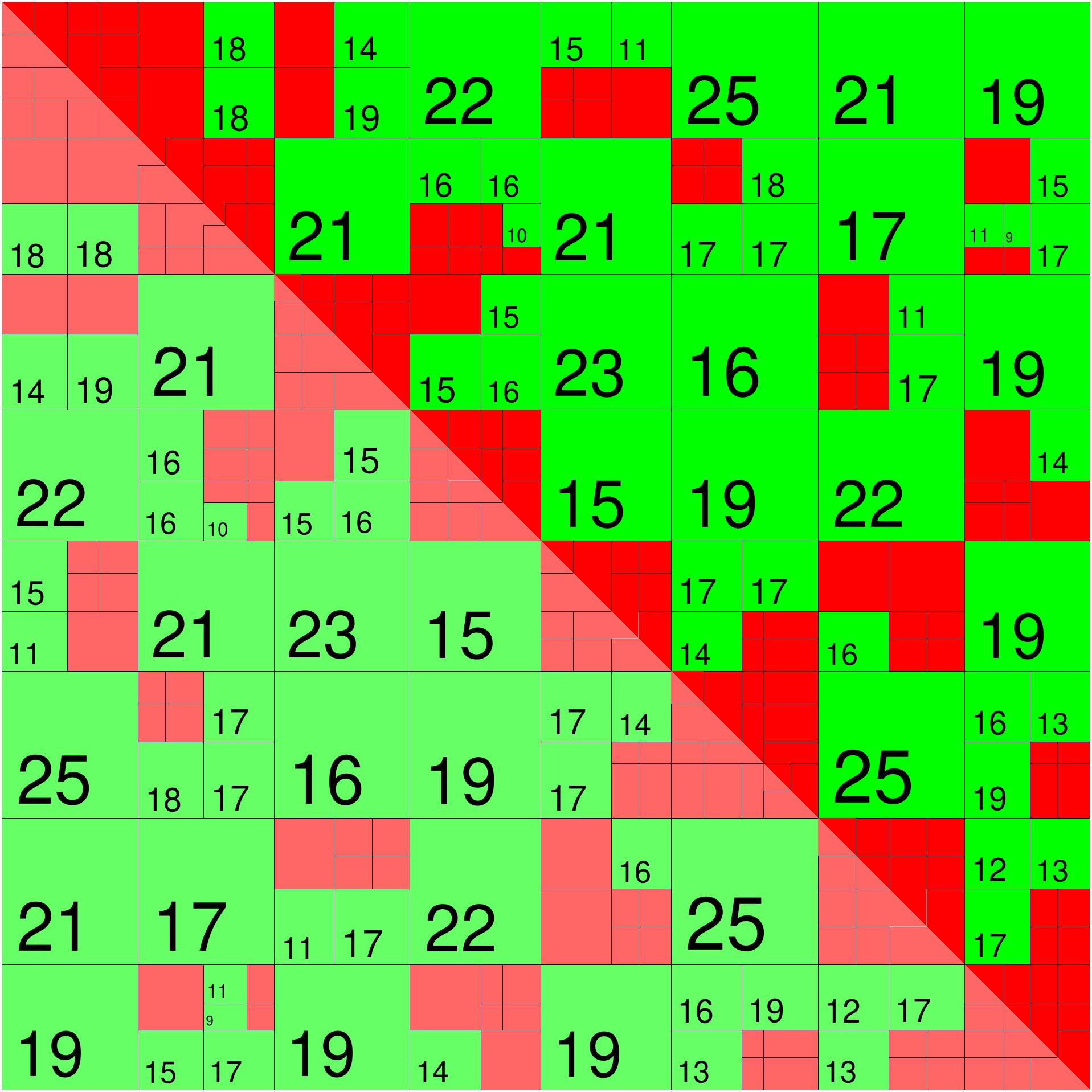}
   \caption{$\eta=1.4,\varepsilon=5e^{-4}$}
  \end{subfigure}
  \begin{subfigure}[b]{0.245\textwidth}
   \includegraphics[width=\textwidth]{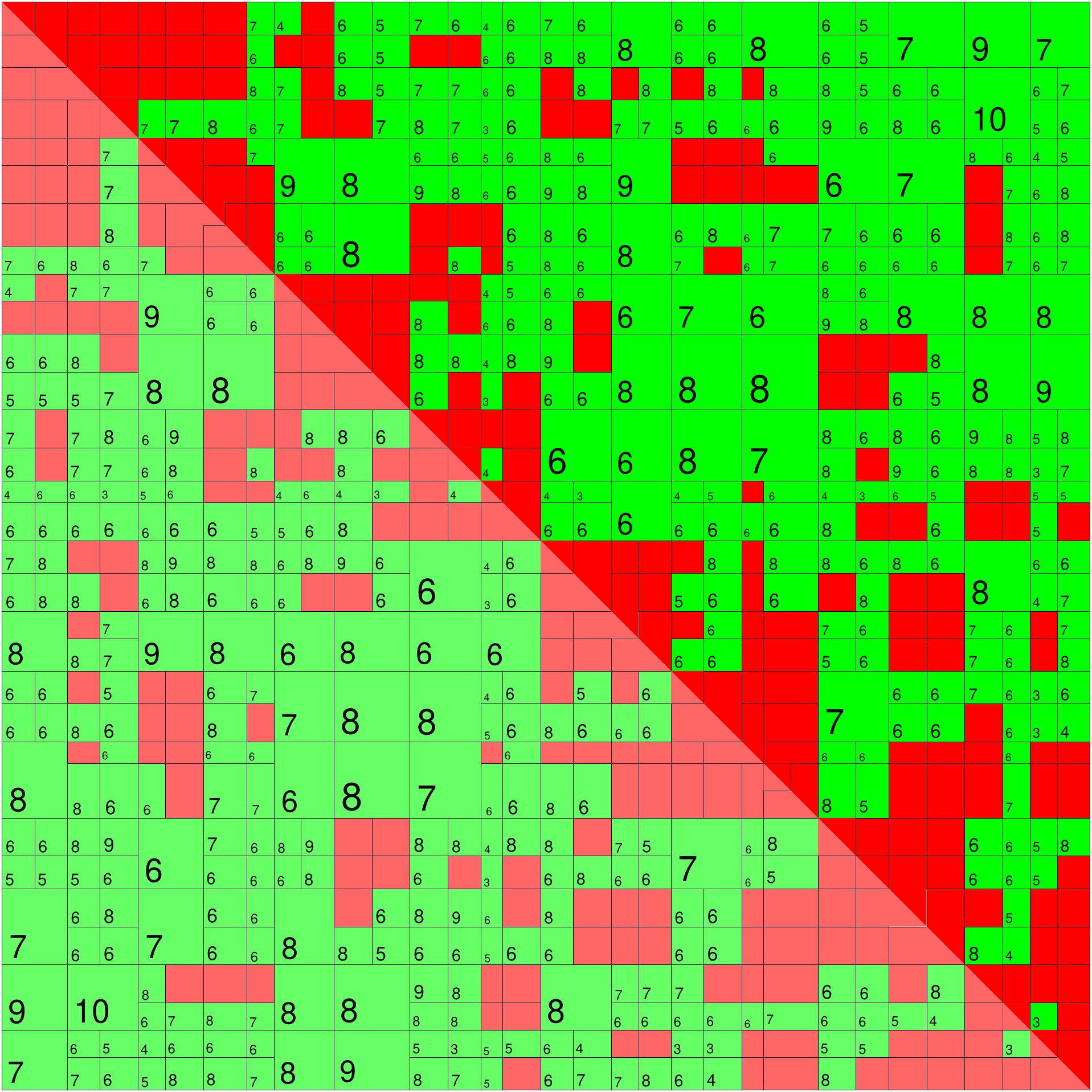}
   \caption{$\eta=0.8,\varepsilon=1e^{-2}$}
  \end{subfigure}
  \begin{subfigure}[b]{0.245\textwidth}
   \includegraphics[width=\textwidth]{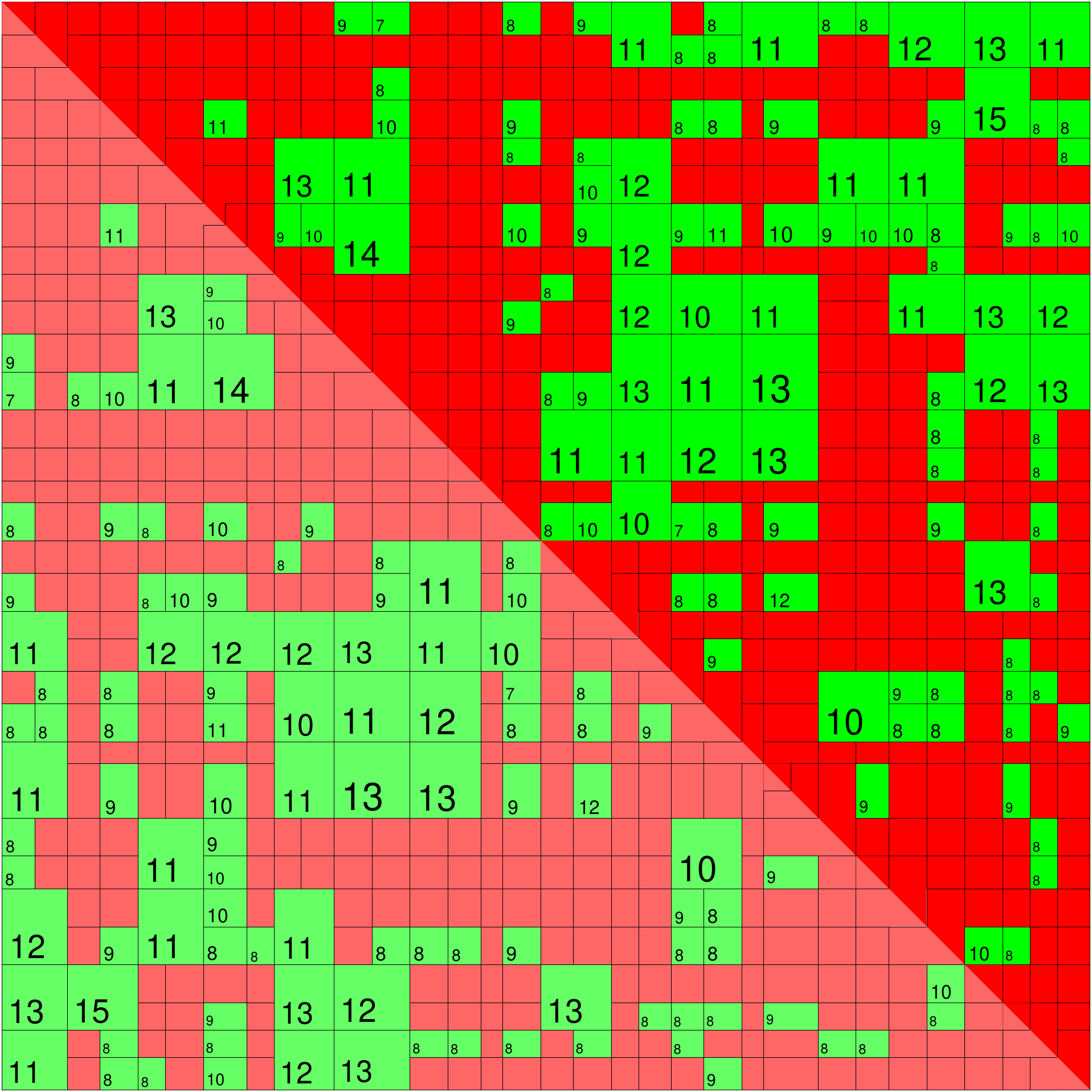}
   \caption{$\eta=0.8,\varepsilon=5e^{-4}$}
  \end{subfigure}
  
  \caption{Representation of the $\mathcal{H}$-matrix approximation of matrix $A$ (200 inclusions) for differents set of parameter $\varepsilon$ and $\eta$. In each admissible block (green) the number corresponds to the rank. Red blocks corresponds to a full storage.}
  \label{fig:hmatrixplot}
\end{figure}

In this paragraph we analyze the influence on the $\mathcal{H}$-matrix approximation of the two parameters:
\begin{itemize} 
\item the admissibility condition parameter $\eta$, see~\eqref{eq:admissibility}, and 
\item the prescribed accuracy for the block approximation $\varepsilon$, see~\eqref{eq:errorA}.
\end{itemize}
The results presented on figure \ref{fig:hmatrixplot} show that $\eta$ influences the block partition: smaller $\eta$ leads to smaller size of the admissible blocks. The choice of $\varepsilon$ influences the rank of the approximation of the blocks: with small $\varepsilon$, the ranks are lower and there are more full-storage blocks (red).

\begin{figure}[h]
  \centering
  \includegraphics[width=0.4\textwidth]{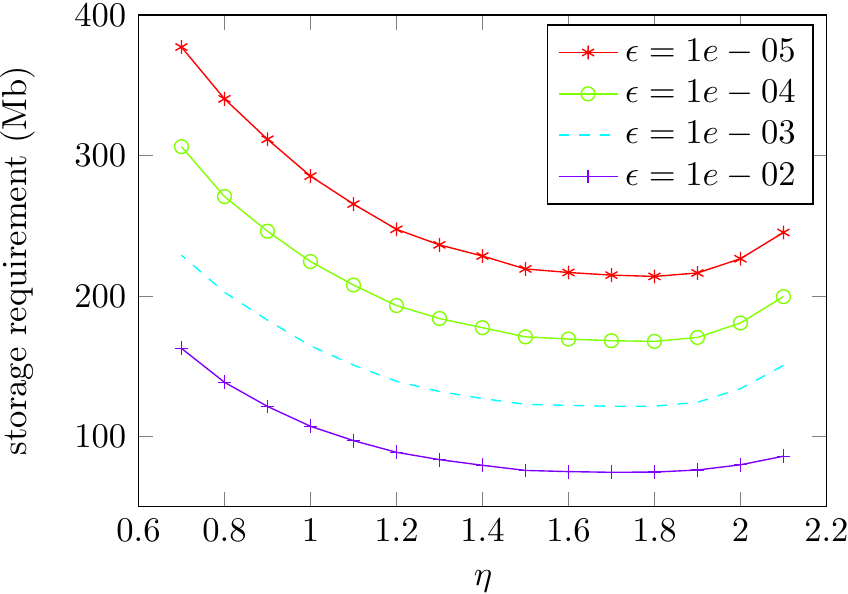}
  \caption{Memory consumption as function of $\eta$ for different precisions $\varepsilon$ (5000 inclusion).}
  \label{fig:memory}
\end{figure}

For a given precision $\varepsilon$, we can not say \textit{a priori} if the compression rate will be better when $\eta$ is large (that corresponds to larger blocks with higher rank) rather than small (which leads to smaller blocks with lower rank).
On figure \ref{fig:memory} we see that $\eta=1.7$ leads to an optimal compression rate, that is a compromise between the size of the blocks and their rank.
\begin{rmrk}
 Note that a good choice choice for $\eta$ depends on the underlying integral equation, and more precisely on the kernel $g$~\cite{Bebendorf2008}.
 Thus $\eta=1.7$ is not necessary the optimal choice for other integral equations.
\end{rmrk}

\subsection{Error analysis}

\begin{figure}[h]
  \centering
  \includegraphics[width=0.45\textwidth]{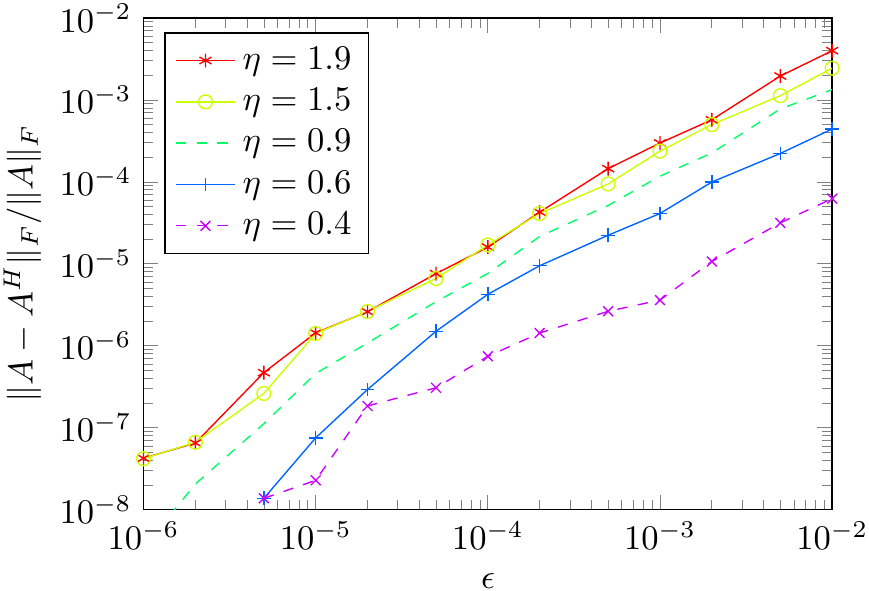}
  \includegraphics[width=0.45\textwidth]{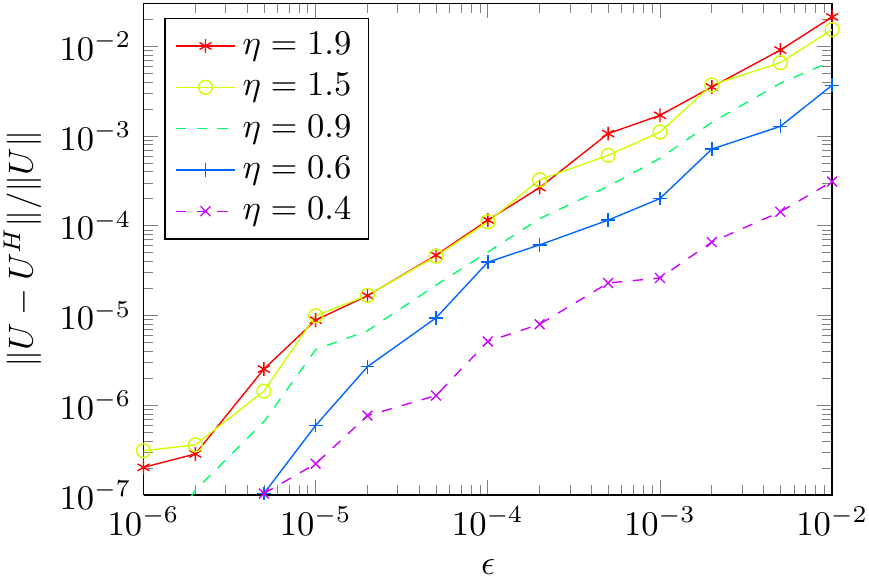}
  \caption{Relative error in the Frobenius norm for the approximation of the operator (left) and relative error on the approximation of the solution (right) as function of $\varepsilon$ (5000 inclusions).}
  \label{fig:erroranalysis}
\end{figure}

We analyze here the error due to the $\mathcal{H}$-matrix approximation.
First we note that condition \eqref{eq:errorA} implies that the approximation $A^H$ has a relative error $\varepsilon$ in the Frobenius norm: $\Vert A-A^H\Vert_F \leq \varepsilon \Vert A\Vert_F$.
Even if the heuristic stopping criterion \eqref{eq:stagnation} of the ACA algorithm with partial pivoting does not ensure \eqref{eq:errorA}, the curves of figure \ref{fig:erroranalysis} shows that in our example, the precision $\varepsilon$ is reached for $A^H$.
Since the Frobenius norm is not an operator norm, we can not directly control the error on the solution by $\varepsilon$.
Denote by $U^H$ the solution of the linear system:
\begin{equation}
 A^H U^H = B. \label{eq:AHUHB}
\end{equation}
We use here a preconditioned conjugate gradient to solve \eqref{eq:AHUHB}. The preconditioner is an incomplete $\mathcal{H}$-LU factorization~\cite{Bebendorf2008}. The iterations are stopped when the residual norm archive a tolerance of $10^{-14}$.
We see on figure \ref{fig:erroranalysis}, using the geometry with 5000 inclusions, that the error on the solution $\Vert U-U^H\Vert/\Vert U \Vert$ is proportional to $\varepsilon$.
This means that if we want an approximation of the solution with respect to a prescribed accuracy, we can impose in practice the same accuracy on the $\mathcal{H}$-matrix approximation of the operator.

\subsection{Scalability}

\begin{figure}[h]
  \centering
  \begin{subfigure}[b]{0.45\textwidth}
   \includegraphics[width=\textwidth]{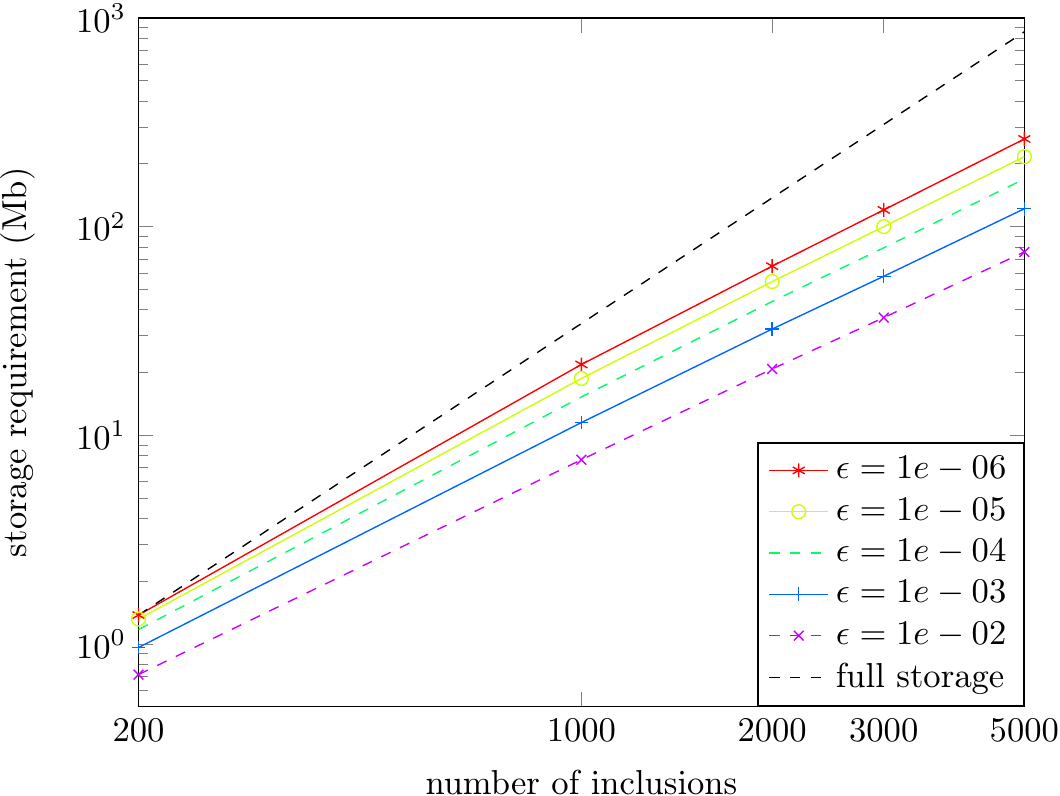}
   \caption{Storage requirement as function of the number of inclusion, for different precision $\varepsilon$ (with $\eta=1.7$)}
   \label{fig:memory_scal}
  \end{subfigure}
  \begin{subfigure}[b]{0.45\textwidth}
   \includegraphics[width=\textwidth]{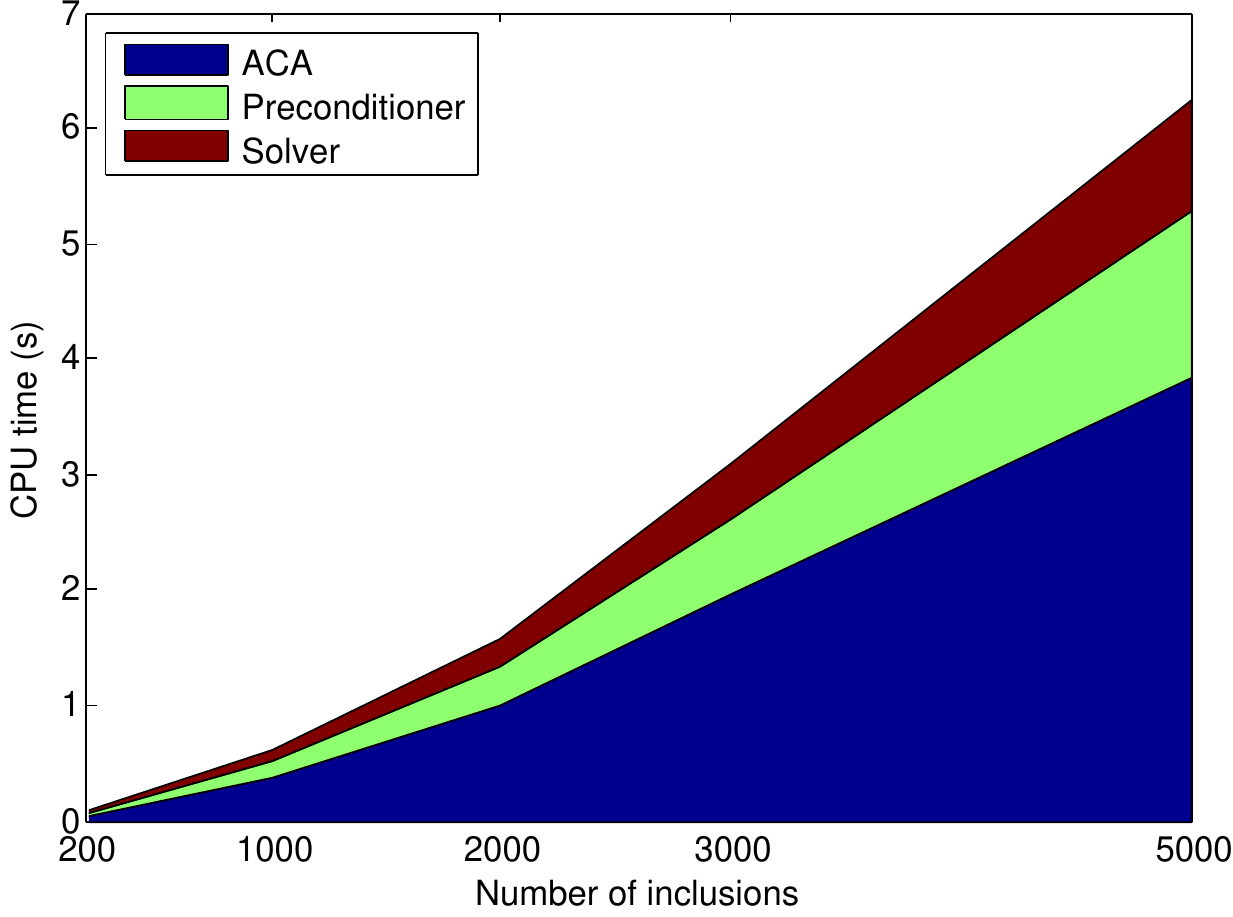}
   \caption{CPU time as function of the number of inclusion (with $\varepsilon=10^{-3}$ and $\eta=1.7$)}
   \label{fig:time_scal}
  \end{subfigure}
  \caption{Evolution of computational time and memory requirement as function of the number of inclusion.}
  \label{fig:scale}
\end{figure}

We assess now the ability of the $\mathcal{H}$-matrix to solve the equivalent inclusion equation for a large number of inclusions.
Figure \ref{fig:memory_scal} shows the memory size for the approximation $A^H$ for different precision $\varepsilon$, with $\eta=1.7$.
We observe that the memory scales as $\mathcal{O}(N^{1.4})$ where $N$ denotes the number of inclusion.
In comparison with the memory scaling as $N^{2}$ which is necessary when using the usual full storage, we see that the $\mathcal{H}$-matrix has a better scaling.
In particular, for $N=5000$, the memory for the $\mathcal{H}$-matrix approximation represent only $8.7\%$ of the memory needed for the full storage.

Furthermore, we present on figure \ref{fig:time_scal} the computational timings on a laptop with a 2.5GHz Intel Centrino 2 for:
\begin{enumerate}
 \item the creation of the $\mathcal{H}$-matrix approximation using ACA with partial pivoting,
 \item the computation of the $\mathcal{H}-$LU preconditioner,
 \item the resolution of the system by a preconditioned conjugate gradient.
\end{enumerate}
We observe that the assembly of the matrix (ACA) is the dominating factor. Also the total timing scales almost linearly with $N$.

\section{Reformulation as a boundary integral problem and the BEM method}
\label{sec:BEM}
We propose now to study a second application of the hierarchical matrix format to solving periodic corrector problems arising in homogenization. This new approach relies on a different formulation of the corrector problem. Rather than using the Lippman--Schwinger equation presented in Section~\ref{sec:equivalentinclusion}, we write a boundary integral problem equivalent to the original corrector problem~\eqref{eq:pb1} with periodic boundary conditions.

\subsection{Motivation}
The equivalent inclusion method presented in the previous sections is interesting by its simplicity: the use of analytic expressions makes for a simple formulation with a variational interpretation as a Galerkin approach which yields a bound for the effective modulus~\cite{Brisard11}. However, the use of one degree of freedom per inclusion is insufficient to obtain a reasonable approximation in materials with densely packed inclusions, which is the case for many applications such as modeling concrete. Indeed, the error resulting from this method is unknown and cannot be directly improved. It is thus necessary to improve on this approach to obtain a better estimate of the effective parameters of the homogenized material.

One possible improvement suggested in~\cite{Moschovidis1975} is to use polynomial approximations of increasing degree inside each inclusion as a Galerkin method. Analytic expressions can still be derived to account for the interactions between inclusions, but the resulting estimate is biased by the use of different boundary conditions implied by the approximation of the Green's function~\eqref{eq:approxgamma}.

We propose to study here as an alternative a rather different approach which is based on the classical reformulation of elliptic boundary value problems as boundary integral equations, which can be solved numerically by the use of boundary element methods (BEM). These methods are known to easily handle jumps in the parameters in a complicated geometry, avoid having to discretize the entire domain in a manner that accurately resolves the geometry of the inclusions, can achieve rapid convergence, and have a robust mathematical foundation~\cite{TextbookBEM, BarnettGreengard}.

 Up to our knowledge, this work is the first attempt at applying the boundary element method to corrector problems arising in periodic or random homogenization. Specific difficulties in this framework involve in particular dealing with:
\begin{itemize}
\item periodic boundary conditions,
\item a three-dimensional setting,
\item complex geometries involving large numbers of randomly distributed inclusions in a representative periodic cell.
\end{itemize}
A closely related problem is the calculation of band structure in periodic crystals or materials, to which boundary integral methods have been recently applied, see e.g.~\cite{Yuan2008}. In particular, our approach is inspired by ideas proposed in~\cite{BarnettGreengard} in a two-dimensional setting.

\subsection{Reformulation of the corrector problem as a boundary integral problem}

\begin{figure}[ht]
\begin{center}
\includegraphics[width=.4\textwidth]{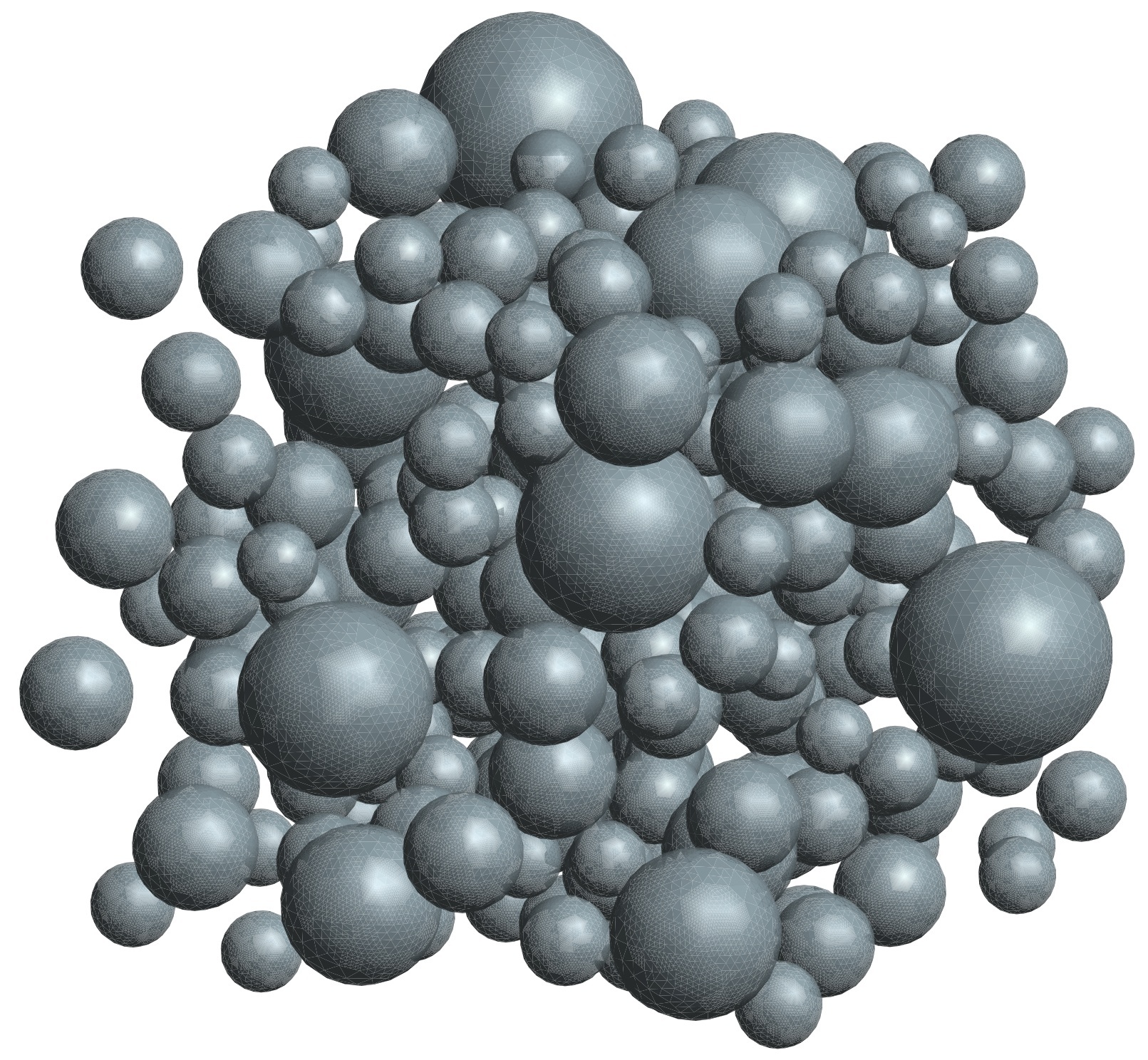}
\caption{Example of the geometry of a set of inclusions.}
\label{fig:geometryinclusions}
\end{center}
\end{figure}

In this section we derive a boundary integral formulation for the corrector problem associated to the Laplace equation. Let $\Omega = (-1/2, 1/2)^d$ be an open unit periodic cell in $\mathbb{R}^d$, $d = 2,3$. We decompose $\Omega$ as a set of closed inclusions $\Omega_{int} \subset \Omega$ and a connected matrix $\Omega_{ext} = \Omega \setminus \overline{\Omega_{int}}$, see figure~\ref{fig:geometryinclusions}. Let $\Gamma = \partial \Omega_{int}$ be the boundary of the set of inclusions.
\begin{rmrk}
It is not necessary to suppose that the inclusions are distributed strictly inside the periodic cell.
\end{rmrk}
Let $\kappa(\mathbf{x})$ be a scalar diffusion coefficient which takes two different values $\kappa_{int}$, $\kappa_{ext} > 0$ respectively in the inclusions and in the matrix:
\begin{equation}\label{def:diffcoefficient}
 \kappa(\mathbf{x}) = \begin{cases}
 \kappa_{int} & \text{ in }  \Omega_{int},\\
 \kappa_{ext} & \text{ in } \Omega_{ext} .
 \end{cases}
\end{equation}
As in Section~\ref{sec:equivalentinclusion}, we are interested in the periodic corrector $u \in H^1(\Omega) / \mathbb{R}$ solution of the following boundary value problem:
 \begin{equation}\label{syst:percorrproblem}
 \left \{ \begin{aligned}
 - \mathrm{div} ( \kappa ( E +  \nabla u )) = 0, &&\text{ in } \Omega,\\
 u \text{ is }\Omega\text{-periodic.}
 \end{aligned} \right.
 \end{equation} 
It is well--known that problem~\eqref{syst:percorrproblem} has a unique solution, see e.g.~\cite{Allaire92}. Note that the solution $u$ is harmonic in $\Omega \setminus \Gamma$ and satisfies continuity conditions:
\begin{equation}\label{eq:matchingconditions}
u,\ \kappa (\nabla u + E) \cdot \mathbf{n}^{ext}\text{ continuous across }\Gamma,
\end{equation}
where $\mathbf{n}^{ext}$ is the exterior normal to $\Gamma$. An elegant approach for the integral representation of periodic fields involves the use of the Green's function satisfying the periodic boundary conditions exactly. Let $\delta_0$ be the Dirac delta function centered at the origin. The periodic Green's function $G_{per}$ is then defined by $-\Delta G_{per} = \delta_0 - 1$ over $\Omega$ (note that the right-hand side must have zero mean over the periodic cell). Let us introduce the single-layer potential~\cite{TextbookBEM} as an operator $H^{-1/2}(\Gamma) \to H^1(\Omega)$ defined by:
\begin{equation}\label{def:singlelayerpotential}
( \Psi_{SL}\sigma)(\mathbf{x}) =  \int_\Gamma  G_{per}(\mathbf{x} - \mathbf{y})  \sigma(
\mathbf{y}) \mathrm{d}s(\mathbf{y}) \quad \text{ for } \mathbf{x} \in \Omega \setminus \Gamma.
\end{equation}
We can represent the solution to the PDE problem~\eqref{syst:percorrproblem} as $u = \Psi_{SL} \sigma$ with a single-layer density $\sigma$. It remains to solve for the density $\sigma \in H^{-1/2}(\Gamma)$ so that the matching conditions~\eqref{eq:matchingconditions} are satisfied. We use the standard jump relations for single layer potential~\cite{TextbookBEM} to write
\begin{equation}\label{eq:BIE1}
\kappa_{ext} \left ( -1/2 Id + K'_0 \right ) \sigma - \kappa_{int} \left ( 1/2 Id + K'_0 \right ) \sigma = -(\kappa_{ext} - \kappa_{int}) E \cdot \mathbf{n}^{ext} \quad \text{ on } \Gamma,
\end{equation}
where $ Id$ is the identity operator and $K'_0$ is the adjoint double layer boundary integral operator $ H^{-1/2}(\Gamma) \to  H^{1/2}(\Gamma)$ formally defined as:
\begin{equation}\label{def:adjdoublelayeroperator}
(K'_0 \sigma)(\mathbf{x}) = \int_\Gamma \sigma(\mathbf{y}) \nabla_\mathbf{x} G_{per}(\mathbf{x} - \mathbf{y}) \cdot \mathbf{n}^{ext}(\mathbf{y}) \mathrm{d}s(\mathbf{y}).
\end{equation}
Finally, we obtain the following boundary integral equation:
\begin{equation}\label{eq:BIE2}
\left ( \frac{\kappa_{ext} + \kappa_{int}}{\kappa_{ext} - \kappa_{int}} Id - K'_0\right ) \sigma = E \cdot \mathbf{n}^{ext},
\end{equation}
where the operator appearing on the left--hand side is a compact perturbation of the identity.
\begin{rmrk}
 Note that this reformulation is quite different from the Lippmann--Schwinger integral equation presented in Section~\ref{sec:equivalentinclusion}: we have not introduced a reference medium and the integrals are defined on the boundary of the inclusions instead of the whole volume. There is also no equivalent to the Hashin--Shtrikman variational principle in this new setting. 
\end{rmrk}
 
\subsection{Solution with the Boundary Element Method and numerical results}
Since the goal of this study is to explore the use of the hierarchical matrix format in the context of the solution of the corrector problem, we do not fully detail here the numerical discretization method and we refer to a forthcoming paper~\cite{futur_BEMHom} for the full description of the approach. We limit ourselves to the case of a smooth boundary $\Gamma$ in dimension $d = 3$.

\subsubsection{Discretization}

To solve numerically the boundary integral equation~\eqref{eq:BIE2}, we propose to use the classical Boundary Element Method. Note that this approach gives rise to full matrices for the representation of the operators, and to reduce the complexity of the computation, compression techniques such as multipole expansions (e.g.~the fast multipole method) or hierarchical matrices can be used; this latter approach is investigated here.

 A Galerkin approach is used to discretize~\eqref{eq:BIE2} given a mesh of the boundary $\Gamma$. Functions in $H^{-1/2}(\Gamma)$ are approximated by constant-by-cell functions. The main difficulty is the computation of values of the periodic Green's function necessary to the discretization of the boundary integral operator $K'_0$ defined by~\eqref{def:adjdoublelayeroperator}. Indeed, while the free--space Green's function is given by an explicit expression, as seen in Section~\ref{sec:equivalentinclusion}, the periodic Green's function is not known explicitely. It is rather given by a Fourier series,
 $$
 G_{per}(\mathbf{x}) = \sum_{\mathbf{k} \in \mathbb{Z}^d;\ \mathbf{k} \neq 0 } \frac{e^{2i\pi \mathbf{k} \cdot \mathbf{x}}}{\vert \mathbf{k} \vert ^2}.
 $$

 Because such sums converge too slowly to be numerically useful, a number of schemes have been devised for their evaluation, see e.g.~\cite{LatticeSums, EwaldSummation}. We use here a different, new approach, based on an idea proposed recently by Barnett and Greengard~\cite{BarnettGreengard} for the Helmholtz two--dimensional case. Observe that the periodic Green's function can be represented in $\Omega$ as
 \begin{equation}\label{eq:GreenFunctionRepresentation}
 G_{per}(\mathbf{x}) = G_\infty(\mathbf{x}) + \frac{\vert \mathbf{x} \vert ^2}{6} +  G^r_{per}(\mathbf{x}),
 \end{equation}
 where $G_\infty$ is the free-space Green's function, introduced in section~\ref{sec:HashShtrik}, and $ G^r_{per}$, the "regular" part of the periodic Green's function, is a smooth, harmonic function in $\Omega$. It can thus be expanded on the basis of the solid spherical harmonics as a series converging uniformly in $\Omega$:
 \begin{equation}\label{eq:GreenFunctionExpansion}
  G^r_{per}(\mathbf{x}) = \sum_{l = 0}^\infty \sum_{m = -l} ^l \alpha_l^m \Phi_l^m(\mathbf{x}),
 \end{equation}
 where $\Phi_l^m$ is the regular solid harmonics of degree $l$ and order $m$ and the $\alpha_l^m$ are scalar coefficients. An approximation to $G^r_{per}$ can then be computed by truncating the expression~\eqref{eq:GreenFunctionExpansion} up to some degree $L \in \mathbb{N}$, given that we tabulate beforehand the coefficients of the expansion $\alpha_l^m$ for $l \leq L$.
 
 \begin{figure}[ht!]
\begin{center}
\includegraphics[width=.45\textwidth]{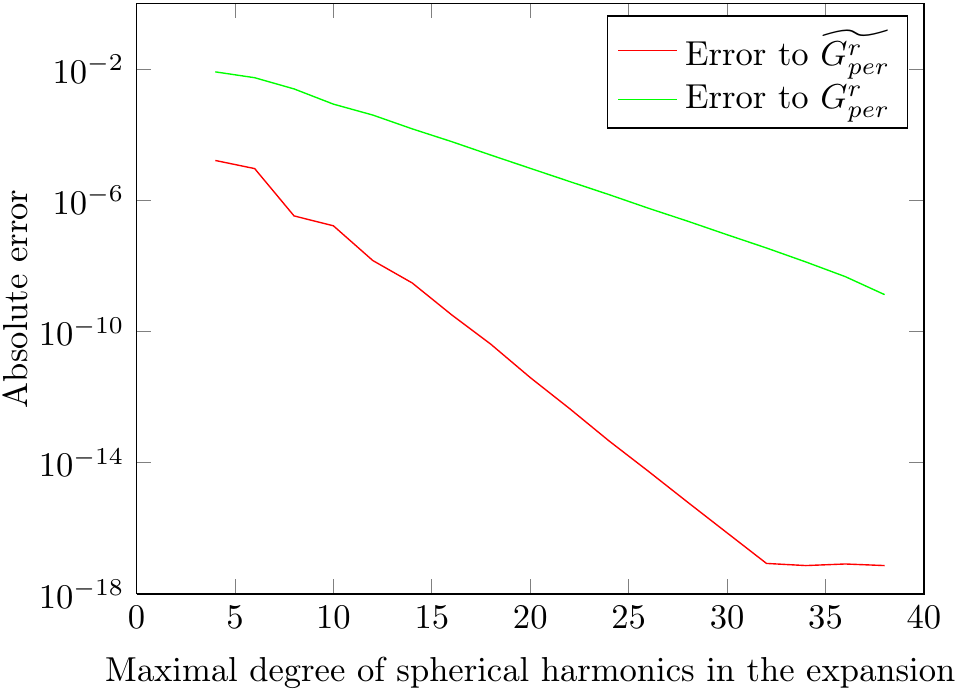}
\caption{Convergence in absolute error, sampled in $\Omega$, for the periodic Green's function scheme.}
\label{fig:convrategreen}
\end{center}
\end{figure}
  To achieve this, we construct a small linear system that enables us to compute an approximation to the coefficients $\alpha_l^m$ easily, as in~\cite{BarnettGreengard}, section 3.2. The general idea is to enforce numerically the periodic boundary conditions on the approximation to $G_{per}$ based on the representation~\eqref{eq:GreenFunctionRepresentation} and ~\eqref{eq:GreenFunctionExpansion}. In practice, we use a least--squares algorithm, minimizing the $L^2$-norm of the deviation from periodicity of the function and its normal derivative on $\partial \Omega$, sampled at Gaussian quadrature points. We refer to the forthcoming paper~\cite{futur_BEMHom} for a more complete description and analysis of this method.

 This method, implemented in Matlab, allows us to achieve an exponential convergence rate as shown by the green curve in figure~\ref{fig:convrategreen}. It can be further accelerated by removing from the regular part $G^r_{per}$ copies of the free-space Green's function centered in the nearest neighbor cells. We then compute and use the coefficients $\beta_l^m$ in the following expansion:
 \begin{equation}\label{eq:GreenFunctionExpansion2}
  \widetilde{G^r_{per}}(\mathbf{x}) = \sum_{l = 0}^\infty \sum_{m = -l} ^l \beta_l^m \Phi_l^m(\mathbf{x}) = G_{per}(\mathbf{x}) -  \sum_{\mathbf{m} \in \{ -1, 0, 1 \}^3 } G_\infty \left (\mathbf{x} + \sum_{i = 1}^3 m_i \mathbf{e}_i \right ) - \frac{\vert \mathbf{x} \vert ^2}{6},
 \end{equation}
 where the $\mathbf{e}_i$ are the vectors of the standard basis of $\mathbb{R}^3$.
 The convergence rate is much improved as seen in figure~\ref{fig:convrategreen}. We use this approximation below, fixing $L = 9$.

\subsubsection{Implementation}
\begin{figure}[ht!]
\begin{subfigure}{\textwidth}
\centering
\includegraphics[width=0.65\textwidth]{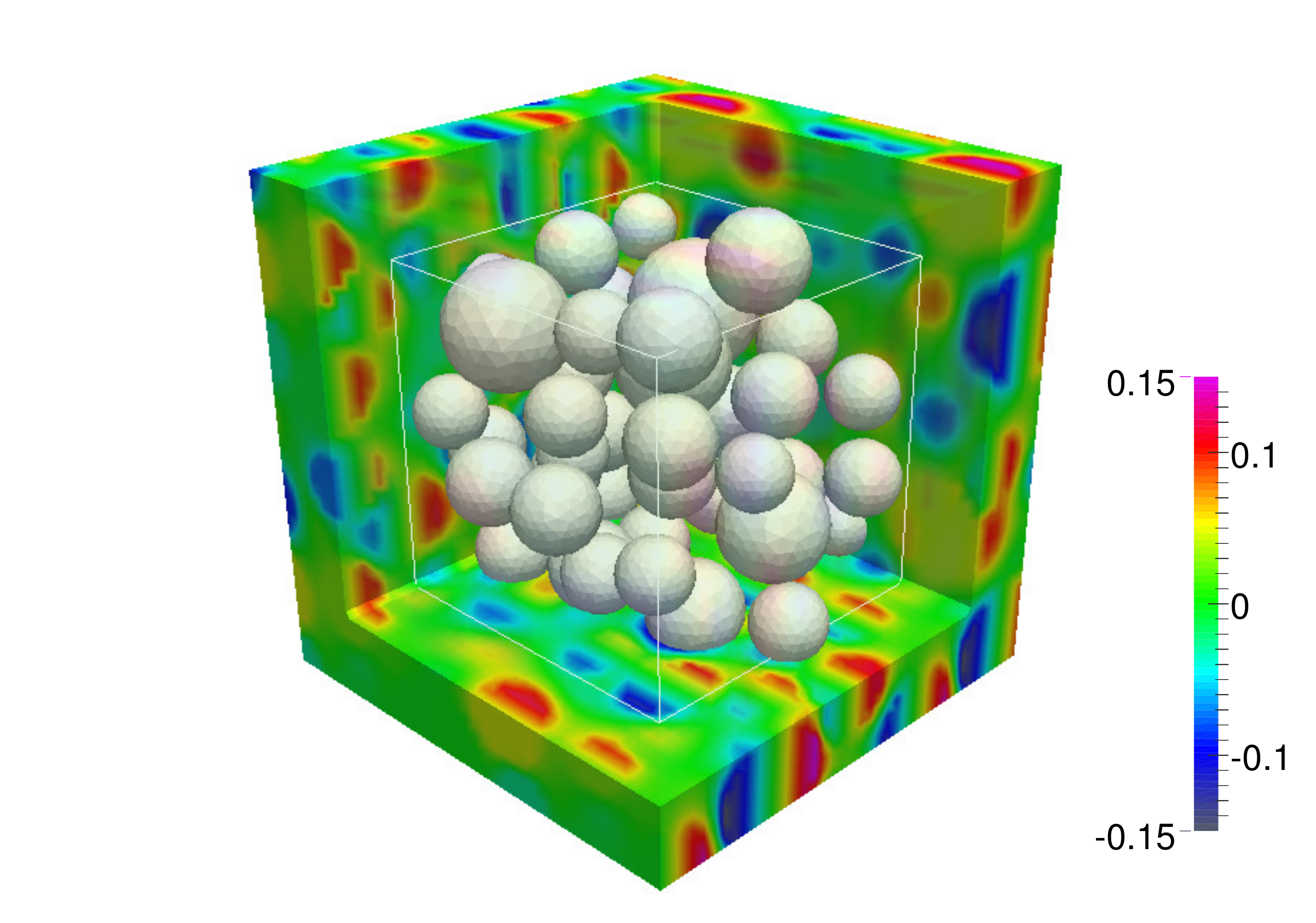}
\caption{Microstructure with 42 inclusions (meshed with 18424 triangles)}
\end{subfigure}
\begin{subfigure}{.49\textwidth}
\centering
\includegraphics[width=\textwidth]{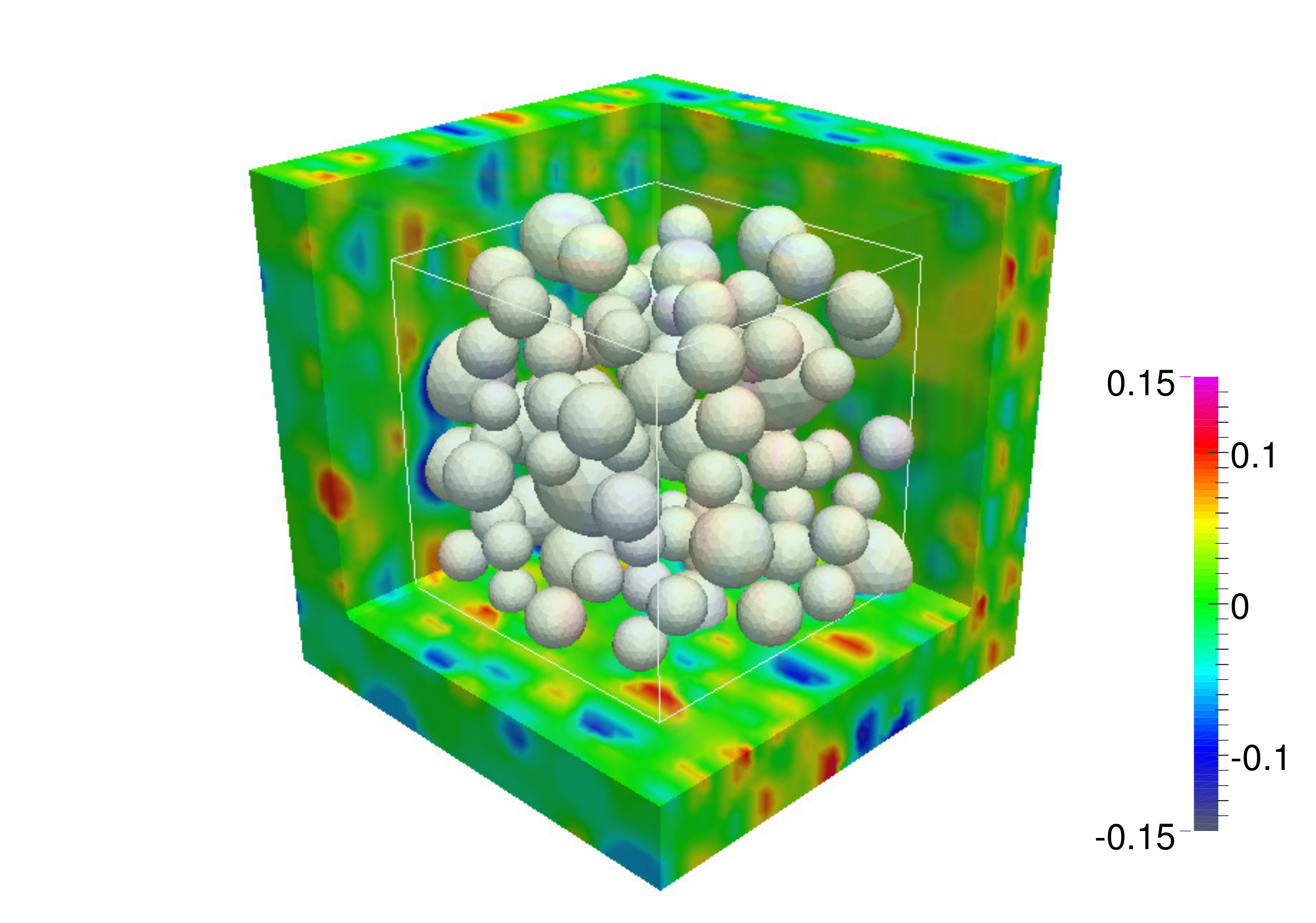}
\caption{Microstructure with 98 inclusions (43034 triangles)}
\end{subfigure}
\begin{subfigure}{.49\textwidth}
\centering
\includegraphics[width=\textwidth]{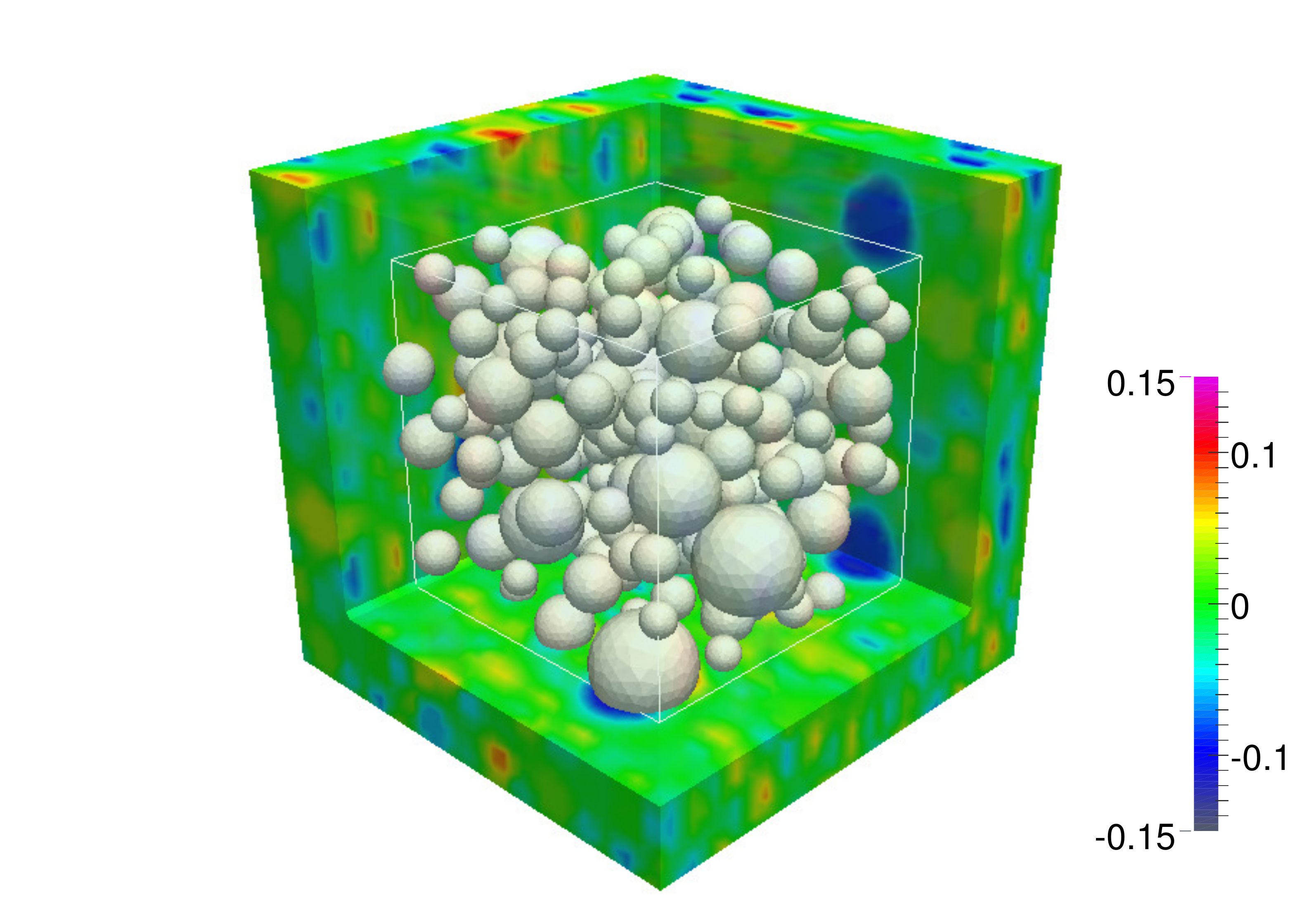}
\caption{Microstructure with 203 inclusion (89140 triangles)}
\end{subfigure}
\caption{Representation of the reconstructed field of the corrector $u$, computed by the boundary element method for the three studied microstructures.}
\label{fig:solutionsBEM}
\end{figure}

We present here some numerical results to illustrate this new method for the resolution of the corrector problem. For the implementation, we used the open-source boundary element library BEM++~\cite{Bempp} interfaced with the library Ahmed~\cite{AHMED} for an efficient representation of the discretized integral operators by $\mathcal{H}$-matrices. All timings are reported for a laptop running the Python / C++ code with a 2.7 GHz Intel Core i7 CPU. Three different microstructures are studied containing respectively 42, 98 and 203 randomly distributed inclusions, with nine different meshes ranging from 1344 to 89140 elements. The diffusion coefficient $\kappa$ takes the value $1$ in the matrix and $100$ in the inclusions. The precision used for the $\mathcal{H}$-matrix compression with ACA algorithm (see section~\ref{sec:hmatrix}) is set to $\varepsilon = 1e-3$.

 Some graphical representations of the solution are shown in figure~\ref{fig:solutionsBEM}. In particular, we can observe the perturbation on the corrector field induced by the periodic copies of the inclusions where the diffusion coefficient is much higher than in the matrix.

\subsection{Analysis of the $\mathcal{H}$-matrix approximation efficiency}
We finally investigate here the scalability of the boundary element approach presented above. As a benchmark, we will also include results obtained with the same data, but replacing the periodic by the free-space Green's function. The behavior and scaling of the method with the number of degrees of freedom is well-known in this case, see e.g.~\cite{Bempp}.

 We present in table~\ref{table:data203} some data collected during the calculations. We observe that the memory storage for the $\mathcal{H}$-matrix representing the integral operator as well as the time used by the solver is slightly increased, but stays of the same order of magnitude in the periodic case in comparison to the benchmark free-space computation. By contrast, the assembly time through the ACA algorithm is much higher in the periodic case. This is due to the much higher cost of evaluation of the periodic Green's function using the representation~\eqref{eq:GreenFunctionRepresentation}. 

\begin{table}[t]
\centering
 \begin{tabular}{|c|c|c|c|c|c|}
 \hline 
mesh size &  Case & storage & compression & ACA time & solver \\ 
  \hline 
44784 triangles &  Free-space & 615 Mb / 15301 MB & 4.02$\%$ & 63s & 76s \\ 
 \hline 
44784 triangles & Periodic &  1012 Mb / 15301 MB & 6.6$\%$ & 3146s & 118s \\ 
 \hline 
89140 triangles &  Free-space & 1374 MB / 60622 MB & 2.26$\%$ & 128s & 190s \\ 
 \hline 
89140 triangles & Periodic &  2210 MB / 60622 MB & 3.6$\%$ & 6402s & 288s \\ 
 \hline 
 \end{tabular} 
\caption{Data for the solution of the boundary element problem in a geometry with 203 inclusions}
\label{table:data203}
\end{table}

\begin{figure}[hb!]
  \centering
  \begin{subfigure}[b]{0.45\textwidth}
   \includegraphics[width=\textwidth]{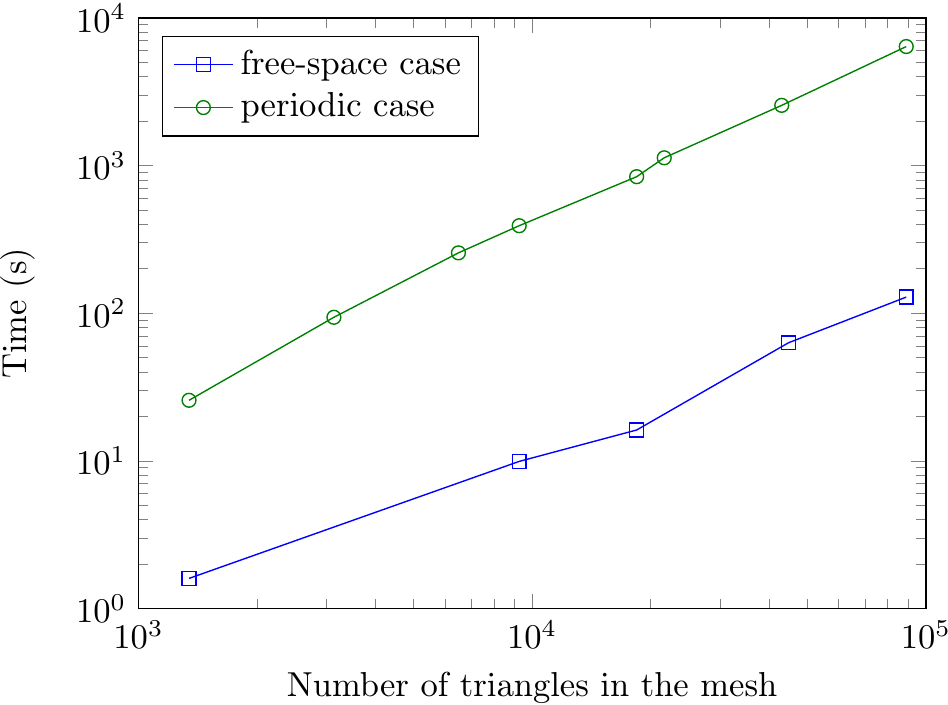}
   \caption{Assembly time for the discrete operator}
   \label{fig:memory_scal_bem}
  \end{subfigure}
  \begin{subfigure}[b]{0.45\textwidth}
   \includegraphics[width=\textwidth]{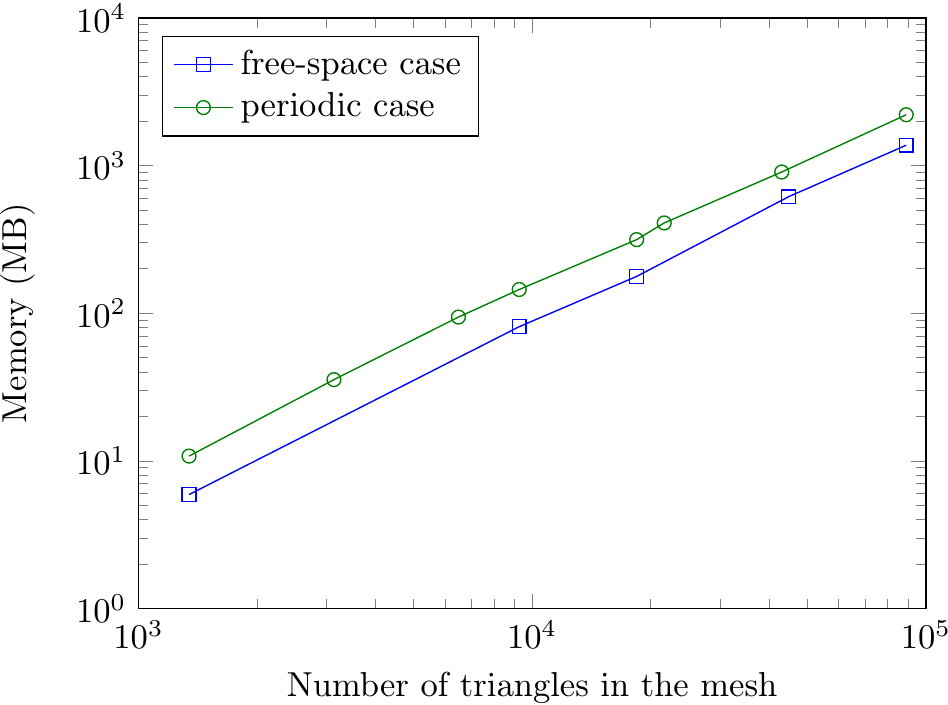}
   \caption{Storage for the discrete operator}
   \label{fig:time_scal_bem}
  \end{subfigure}
  \caption{Computational cost scaling as a function of mesh size for the periodic and the free-space cases.}
  \label{fig:datacompBEM}
\end{figure}

Figure~\ref{fig:datacompBEM} illustrates however that the scaling in performance with problem size is independent of the kernel and also of the number of inclusions in the computational geometry: in all cases we observe a memory and time cost scaling approximately as $\mathcal{O}(N^{1.3})$.

\section{Conclusion}
The $\mathcal{H}$-matrix technique has been successfully applied for the resolution of integral equations arising from corrector problem of thermal homogenization.
The gain of memory and of computational costs provided by the $\mathcal{H}$-matrix format has allowed to deal with a very large number of degrees of freedom in the simulation.
With respect to the equivalent inclusion method, we have shown the feasibility of solving the problem with large numbers of inclusions. However, a more refined discretization of the polarization field is necessary in order to provide accurate predictions of the effective properties, rather than further adding inclusions in the RVE. We refer to~\cite{Brisard11} for details in this direction, for which the $\mathcal{H}$-matrix approach should remain very effective.

We have also proposed a new approach to the solution of the corrector problem by the use of boundary integral equations. The first results are very promising, and we hope to further develop this method which will be presented in details and compared to existing approaches in~\cite{futur_BEMHom}.

\bibliographystyle{plain}
\bibliography{biblio}

\end{document}